\documentclass[11pt]{article}

\usepackage[alphabetic]{amsrefs}
\usepackage{amsmath,amssymb,amsthm,enumerate}
\usepackage[all,cmtip]{xy}
\usepackage[applemac]{inputenc}

\def \a{{\mathfrak a}}

\def \al{\alpha}

\def \be{\beta}

\def \C{{\mathbb C}}

\def \CC{{\cal C}}
\def \CM{\mathcal{M}}
\def \CN{\mathcal{N}}
\def \CO{{\cal O}}
\def \CR{{\cal R}}
\def \CZ{{\cal Z}}
\def \CSch{{\rm CSch}}
\def \df{\ \stackrel{\mbox{\rm\tiny def}}{=}\ }
\def \e{\emph}
\def \Ess{\operatorname{Ess}}
\def \F{{\mathbb F}}
\def \Ga{\Gamma}
\def \ga{\gamma}
\def \GL{\operatorname{GL}}
\def \Hom{\operatorname{Hom}}
\def \Id{\operatorname{Id}}
\def \Im{\operatorname{Im}}

\def \m{{\mathfrak m}}

\def \Mspec{\operatorname{Mspec}}
\def \N{{\mathbb N}}
\def \Nil{\operatorname{Nil}}
\def \O{{\rm O}}
\def \OP{{\rm OP}}
\def \p{{\mathfrak p}}
\def \Per{\operatorname{Per}}
\def \ph{\varphi}
\def \PlogSch{\operatorname{PlogSch}}
\def \q{{\mathfrak q}}
\def \qLR{\quad\Leftrightarrow\quad}
\def \Q{{\mathbb Q}}
\def \q{{\mathfrak q}}
\def \qR{\quad\Rightarrow\quad}
\def \R{{\mathbb R}}
\def \Re{\operatorname{Re}}
\def \red{{\rm red}}
\def \res{\operatorname{res}}
\def \RINGS{{\rm RINGS}}
\def \SES{{\rm SES}}
\def \SL{{\rm SL}}
\def \Sp{{\rm Sp}}
\def \sm{\raisebox{1pt}{$\smallsetminus$}}
\def \spec{\operatorname{spec}}

\def \what{\widehat}
\def \Z{{\mathbb Z}}
\def \ZSch{{\Z{\rm Sch}}}
\def \({\left(}
\def \){\right)}

\newcommand{\tto}[1]
{\stackrel{#1}{\longrightarrow}}

\renewcommand{\sp}
[1]{\left\langle #1\right\rangle}

\newcommand{\ol}
[1]{\overline{#1}}

\newcommand{\stack}
[2]{\genfrac{}{}{0pt}{1}{#1}{#2}}

\newcommand{\cupdot}
{\ensuremath{\mathaccent\cdot\cup}}

\newcommand{\smat}
[4]{\(\begin{smallmatrix}#1 & #2 \\ #3 & #4\end{smallmatrix}\)}

\newtheorem{theorem}{Theorem}[subsection]
\newtheorem{conjecture}[theorem]{Conjecture}
\newtheorem{lemma}[theorem]{Lemma}
\newtheorem{corollary}[theorem]{Corollary}
\newtheorem{proposition}[theorem]{Proposition}
\newtheorem{exmple}[theorem]{Example}
\newenvironment{example}[0]{\begin{exmple}\rm}
{\end{exmple}}
\newtheorem{exmples}[theorem]{Examples}
\newenvironment{examples}{\begin{exmples}\nopagebreak\begin{itemize}\nopagebreak\rm}{\end{itemize}\end{exmples}}
\newtheorem{defi}[theorem]{Definition}
\newenvironment{definition}[0]{\begin{defi}\rm}
{\end{defi}}

\begin{document}

\pagestyle{myheadings} \markright{CONGRUENCE SCHEMES}

\title{Congruence schemes}
\author{Anton Deitmar\\ \ \\
International Journal of Math. Vol. 24, Issue 2 [46 pages] (2013)}
\date{}
\maketitle

{\bf Abstract:}
A new category of algebro-geometric objects is defined which contains the category of monoid schemes, or schemes over $\F_1$, and the category of Grothendieck schemes as subcategories.

$$ $$

%\newpage
\section*{Introduction}

In \cite{KOW}, the authors suggested the idea of a non-additive geometry using the sets of congruences on monoids, where
a congruence is a structure-preserving equivalence relation.
In a case of a ring, all congruences come from ideals. For monoids, however, there are far more congruences than ideals.
Congruences on monoids do not lend themselves well for forming structure sheaves.
Vladimir Berkovich \cites{Berk, Berk2} gives a line of attack to this problem which is influenced by the theory of $p$-adic spaces and designed to give an abstract definition of Berkovich skeletons.

In the present paper we give a construction of structure sheaves on the Zariski site, which is adapted to number theoretical problems arising from the comparison of number fields and function fields.
This construction is a vast generalization of existing $\F_1$-theories, as it contains the the theory of monoid schemes \cite{F1} on the one end and classical algebraic theory, e.g. Grothendieck schemes, on the the other.
It also gives a handy description of Berkovich subdomains and thus contains Berkovich's approach to abstract skeletons.
Further it complements the theory of monoid schemes \cite{F1} in view of number theoretic applications as congruence schemes encode number theoretical information as opposed to combinatorial data which are seen by monoid schemes.

After this paper was (pre-)published, the paper \cite{blueprints} was put on the server.
In the latter paper, ideals are considered instead of congruences. In section \ref{relation} we give a comparison of the two approaches.
I thank Walter Gubler for bringing \cite{Berk} to my attention and many discussions on the subject of this paper. I thank Oliver Lorscheid and Philipp Vollmer for comments on the first draft.

{\small
\baselineskip 0pt
\tableofcontents
}

%: 
\section{Sesquiads}
\subsection{Definition}
In this paper, a \e{ring} will always be commutative with 1 and a \e{monoid} will always be commutative, so it is a set $A$ with an associative and commutative composition $(a,b)\mapsto ab$ and a unit element $1\in A$ with $1a=a$ for every $a\in A$.
A \e{monoid morphism} $\ph:A\to B$ is a map between monoids such that $\ph(aa')=\ph(a)\ph(a')$ for all $a,a'\in A$ and $\ph(1)=1$.
For a monoid $A$, let $A^\times$ denote its \e{unit group}, i.e., the group of all invertible elements.

A \e{zero element} of a monoid $A$ is an element $a_0$ such that $a_0a=a_0$ for every $a\in A$.
If it exists, it is uniquely determined and we write it as $0\in A$.
To a given monoid $A$ we can attach a zero element, defining $A_0=A\cup\{ 0\}$ with the obvious monoid structure.

\begin{center}
{\it For simplicity of presentation, in this paper all monoids have a zero.\\
We further insist that monoid morphisms preserve zero.}
\end{center}

For a monoid $A$ we define the \e{monoidal ring} $\Z A$ as the free $\Z$-module generated by $A$ with multiplication given by the product in $A$.
We further define the \e{reduced monoidal ring}
$$
\Z_0A=\Z A/\Z 0_A,
$$
where $0_A$ is the zero element of $A$.
Note as a special case that if $A=A_1\cup\{ 0\}$ for some monoid $A_1$, then $\Z_0A$ is isomorphic to the monoidal ring $\Z A_1$.

\begin{definition}
An \e{addition} or a 
\e{$+$-structure} on a monoid $A$ is a family $(D_k,\Sigma_k)_{k\in\Z^n, n\ge 2}$ where $D_k\subset A^n$ 
and  $\Sigma_k:D_k\to A$ is a map
such that there exists an injective morphism $\ph:A\hookrightarrow R$ to the multiplicative monoid of a ring $R$ which satisfies $\ph(0)=0$ and $\ph(\sum_k(a))=\sum_{j=1}^nk_j\ph(a_j)$ for every $a\in D_k$.
We further insist that $D_k$ be maximal in the following sense
$$
D_k=\left\{ a\in A^n: \sum_{j=1}^nk_j\ph(a_j)\in \Im(\ph)\right\}.
$$

This implies that the addition is associative and distributive when defined and respects zero, i.e., $a+0=a$ holds for every $a\in A$.
It further implies that addition is cancellative, i.e., if $a+c=b+c$, then $b=a$.
\end{definition}

The maximality of $D_k$ implies that $D_k$ is uniquely determined by $\ph$ and so we can give an alternative definition of a $+$-structure as an equivalence class of injections of a multiplicative monoid $\ph:A\hookrightarrow R$ into rings $R$, where two such embeddings are equivalent, if they define the same addition on $A$. See Section \ref{alternative} for more on this idea.

\begin{definition}
A monoid $A$ together with an addition $(D,\Sigma)$ is called a \e{sesquiad}.
By a \e{sesquiad morphism} $\ph$ we mean a morphism of monoids $A\to B$ such that $\ph\(\sum_k(a)\)=\sum_k(\ph(a))$, or, in a different way of writing,  
$$
\ph\(\sum_{j=1}^nk_ja_j\)=\sum_{j=1}^nk_j\ph(a_j)
$$ 
for all $n\in\N$, $k\in\Z^n$, $a\in D_k$.
\end{definition}

\begin{example}
On the one extreme, a ring is a sesquiad.
On the other, a given monoid $A$ the inclusion $A\hookrightarrow\Z_0A$ defines an addition on $A$, where one only can add zero to any element.
We call this the \e{trivial addition}.
Every monoid morphism becomes a sesquiad-morphism with this addition.
Hence the category of sesquiads contains the categories of rings and monoids as full subcategories.
\end{example}

The name \e{sesquiad} comes from the latin word \e{sesquialter} for ``one and a half''.
This is because a sesquiad has only half an addition, so it has one and a half operations.

\subsection{Universal ring}
The ring representing an addition is not uniquely defined, but there is a universal one, as in the following proposition.

\begin{proposition}\label{prop1.2.1}
To every sesquiad $(A,D,+)$ there is a universal morphism $\ph_A:A\hookrightarrow R_A$ to a ring $R_A$, such that every sesquiad-morphism $\ph:A\to R$ to a ring $R$ factors uniquely over $\ph_A$, i.e., $\ph_A$ induces a functorial isomorphism
$$
\Hom_{\SES}(A,R)\cong\Hom_{\RINGS}(R_A,R),
$$
where the homomorphisms are taken in the category of sesquiads and rings respectively.
The morphism $\ph_A:A\to R_A$ represents the addition.
The ring $R_A$ is uniquely determined up to isomorphism. The morphism $\ph_A$ is uniquely determined up to unique isomorphism.
\end{proposition}

\begin{proof}
Let $(D,\Sigma)$ be an addition on the monoid $A$.
Let $I(A)\subset\Z A$ be the ideal generated by all elements of the form $\sum_{j=1}^nk_ja_j-\Sigma_k(a)$ with $a\in D_k$.
Let $R_A=\Z A/I(A)$ and $\ph_A:A\to R_A$ the canonical map.
For any morphism of sesquiads $\ph:A\to R$, there exists a unique ring homomorphism $\Z A\to R$ factoring $\ph$.
Since $\ph$ is a sesquiad morphism, the morphism $\Z A\to R$ factorizes uniquely over $R_A$.
The injectivity of $\ph_A$ follows a fortiori by applying the above to an injection $A\hookrightarrow R$, which defines the addition on $A$.

Finally, one has to show that the domain $D_k$ is indeed given by the map $\ph_A:A\to R_A$.
For this let $D_k'$ be the set of all $a\in A^n$ such that $\sum_{j=1}^nk_j\ph_A(a_j)\in\Im(\ph_A)$.
We have to show $D_k=D_k'$.
So let $a\in D_k$.
Then $\sum_{j=1}^nk_ja_j-\Sigma_k(a)
\in I(A)$, so $\sum_{j=1}^nk_ja_j\in \ph_A(A)$, i.e., $a\in D_k'$.
For the converse, let $a\in D_k'$.
Let $\ph:A\to R$ be defining the addition on $A$.
Let $\psi$ denote the unique morphism with $\ph=\psi\circ\ph_A$.
Now
\begin{align*}
\sum_{j=1}^nk_j\ph(a_j)&=
\sum_{j=1}^nk_j\psi(\ph_A((a_j))\\
&=
\psi\(\underbrace{\sum_{j=1}^nk_j\ph_A((a_j)}_{\in\Im(\ph_A)}\)\in\Im(\psi\circ\ph_A)=\Im(\ph),
\end{align*}
i.e., $a\in D_k$.
\end{proof}

Let $S\subset A$ be a submonoid, not necessarily containing zero.
The localization $S^{-1}A$ is defined as in the ring case \cite{F1}.
An addition on a monoid $A$ naturally induces an addition on any localization $S^{-1}A$ given by the inclusion $S^{-1}A\hookrightarrow S^{-1}R_A$.
If $S$ contains zero, then $S^{-1} A=0=S^{-1} R_A$, so the interesting case is when $S$ does not contain zero.

\begin{proposition}
The universal ring of the localization $S^{-1}A$ is the ring $S^{-1}R_A$.
The sesquiad $S^{-1}A$ has the following universal property: Every morphism of sesquiads $\phi:A\to B$ which maps $S$ into the unit group $B^\times$, factors uniquely over the localization $S^{-1} A$.
\end{proposition}

\begin{proof}
Let $\ph:S^{-1}A\to R$ be a monoid morphism to a ring $R$.
Via $A\to S^{-1} A\to R$ one gets a ring homomorphism $R_A\to R$, which maps $S$ to the unit group $R^\times$, therefore it induces a ring homomorphism $S^{-1}R_A\to R$.
The uniqueness and 
the universal property of $S^{-1} A$ are clear.
\end{proof}

\begin{lemma}\label{1.2.3}
Let $A$ be a sesquiad and let $B\subset R_A$ another sesquiad such that the inclusion maps
$$
A\hookrightarrow B\hookrightarrow  R_A
$$
are sesquiad morphisms.
Then $R_B=R_A$.
\end{lemma}

\begin{proof}
The inclusion map $A\to B$ induces a ring homomorphism $R_A\to R_{B}$.
The sesquiad morphism $B\hookrightarrow R_A$ induces a ring morphism $R_{B}\to R_A$.
These two morphisms are inverse to each other.
\end{proof}

\begin{proposition}
The map $\rho:A\mapsto R_A$ defines a faithful functor from the category of sesquiads to the category of rings.
This functor preserves surjectivity, but not injectivity
\end{proposition}

\begin{proof}
In the proof of Proposition \ref{prop1.2.1} we actually constructed a universal ring $R_A$ for a given sesquiad $A$.
This means that we have a canonical element within the isomorphy class of all universal rings.
With this convention the map $\rho:A\to R_A$ is well-defined without making any choices.
Now let $\ph:A\to B$ be a morphism of sesquiads.
The defining property of the universal ring $R_A$ gives a unique morphism $\ph_R:R_A\to R_B$ of rings making the diagram
$$
\xymatrix{
A\ar[r]^\ph\ar[d]& B\ar[d]\\
R_A\ar[r]^{\ph_R}& R_B
}
$$
commutative.
This means that $A\mapsto R_A$ defines a faithful functor.
Let $\ph:A\to B$ be a surjective morphism of sesquiads, then the induced ring homomorphism $R_A\to R_B$ is surjective since the ring $R_B$ is generated by the subset $B$.
A counterexample for injectivity is provided by a finite field $\F$.
let $A$ be the sesquiad you get from $\F$ when you take the multiplication  from $\F$ but equip $A$ with the trivial addition.
The identity map $\ph:A\to\F$ is a morphism of sesquiads, but $R_A$ is an infinite ring, whereas $R_\F=\F$.
\end{proof}

\begin{definition}
An injective sesquiad morphism $\ph:A\hookrightarrow B$ is called an \e{embedding}, if it is an isomorphism to a subsesquiad of $B$, i.e., if the addition on $A$ is defined via $\ph$.
Writing $D_{k,A}$ and $D_{k,B}$ for the domains of definition for the additions on $A$ and $B$ respectively, and considering $A$ as a subset of $B$, this property is equivalent to the inclusion being a sesquiad morphism and
$$
D_{k,A}=D_{k,B}\cap A^n\quad\text{for all  }n\in\N,\ k\in\Z^n.
$$
\end{definition}

\begin{proposition}\label{prop1.2.7}
\begin{enumerate}[\rm (a)]
\item The composition of embeddings is an embedding.
If a composition $\ga\circ\al$ of sesquiad morphisms is an embedding, then $\al$ is an embedding.
Consequently, if in the commutative diagram
$$
\xymatrix{
A \ar[r]^\al \ar[d]_\beta &B \ar[d]^\ga\\
C\ar[r]^\delta & D
}
$$
the arrows $\be,\ga,\delta$ are embeddings, then so is $\al$.
\item A morphism of sesquiads $\ph:A\to B$ is an embedding if and only if the induced ring homomorphism $\ph_R:R_A\to R_B$ is injective.
\item Let $A$ be a sesquiad and $R$ a ring.
Suppose that $\ph:A\hookrightarrow R$ defines the addition on $A$.
Then the induced map $R_A\to R$ is injective.
\item If $\ph:A\hookrightarrow R$ defines the addition of $A$ and the ring $R$ is generated by $\ph(A)$, then $R\cong R_A$.
\end{enumerate}
\end{proposition}

\begin{proof}
(a) is easy. For (b) assume $\ph$ is an embedding,
let $\al\in\ker\ph_R$ and write $\al$ as
$$
\al=\sum_{j=1}^nk_ja_j
$$
for some $a_j\in A$ and $k_j\in\Z$.
This means that $\sum_{j=1}^nk_j\ph(a_j)=0$, but as $\ph$ is an embedding, this means that $\Sigma_k(a)$ is defined and is equal to zero, so $\al=0$.
For the converse direction assume that $\ph_R$ is injective. 
Then the addition on $A$ is defined via $A\hookrightarrow R_B$, as is the addition on $B$, so $\ph$ is an embedding.
Finally, (c) is only a special case of (b) and (d) is an easy consequence of (c).
\end{proof}

\begin{example}
In Arakelov theory one faces the difficulty that at the place $\infty$ there is no proper counterpart of the ring on integers $\Z_p$ at the finite place $p$.
For many purposes the set $\{x\in\R:|x|\le 1\}=[-1,1]$ will do, which is why it is often denoted as $\Z_\infty$.
Then $\Z_\infty$ is not  ring, but it is a sesquiad with the addition of the ambient field $\R$ of real numbers.
One naturally speculates, if $\R$ might indeed be its universal ring, and this is the case as we show next.
\end{example}

\begin{corollary}
The universal ring of the sesquiad $\Z_\infty$ is $\R$.
\end{corollary}

\begin{proof}
This follows from Proposition \ref{prop1.2.7}, part (d), as $\R$ is generated as a ring by $\Z_\infty$.
\end{proof}

\subsection{Ideals and congruences}
Let $A$ be a sesquiad.
For a subset $S\subset A$ we say that $S$ is \e{closed under addition}, if for every $n\ge 2$ and $k\in\Z^n$, $s\in S^n$ such that $\sum_k(s)$ is defined in $A$, we have $\sum_k(s)\in S$.
Note that a subset $S$ which is non-empty and closed under addition, always contains zero, as $0\cdot s+0\cdot s =0$ is defined in $A$.

A subset $\a$ of $A$ is called a (sesquiad-) \e{ideal}, if there exists an ideal $I$ of the ring $R_A$ with
$$
\a= I\cap A.
$$
If $\a$ is an ideal of the sesquiad $A$, we can choose $I$ to be the ideal $(\a)$ generated by $\a$.

\begin{lemma}
A non-empty subset $\a$ of a sesquiad $A$ is an ideal if and only if 
\begin{itemize}
\item $\a$ is closed under addition, and
\item $\a A\subset\a$.
\end{itemize}
\end{lemma}

\begin{proof}
For a given ideal $\a$ the two properties are obvious.
Conversely, let $\a$ be a subset of $A$ having the two properties.
Let $I$ be the ideal of the ring $R_A$ generated by $\a$.
We claim that $\a= I\cap A$.
The inclusion ``$\subset$'' is clearly satisfied. 
For the other direction, let $\al\in I\cap A$.
Then $\al=\sum_{j=1}^n r_ja_j$ for some $r_j\in R_A$ and $a_j\in \a$.
Any $r_j$ can be written as a sum $\sum_{\nu=1}^{m_j}k_{j,\nu} a_{j,\nu}$ with $k_\nu\in\Z$ and $a_{i,\nu}\in A$.
Plugging this sum in and noting that $a_{j,\nu}a_j\in\a$ we find that $\al$ can be written as $\al=\sum_{j=1}^mk_j,\al_j$ for some $k_j\in\Z$ and $\al_j\in\a$.
As $\a$ is closed under addition, we conclude $\al\in\a$.
\end{proof}

An ideal $\a$ is called \e{prime}, if its complement $A\sm\a$ is closed under multiplication.
Let $\spec_z A$ be the \e{Zariski spectrum} of $A$, which is the set of all prime ideals with the topology generated by all sets of the form
$$
D(f)=\{\p\in\spec_z A: f\notin\p\}.
$$

\begin{definition}
A \e{congruence} on a sesquiad $A$ is an equivalence relation $\CC\subset A\times A$ such that there is a sesquiad structure on $A/\CC$ making the projection $A\to A/\CC$ a morphism of sesquiads.
This  condition implies $x\sim_\CC y\Rightarrow xz\sim_\CC yz$ for all $x,y,z\in A$ and $x+z\sim_\CC y+z$ if $x+z$ and $y+z$ are defined.
\end{definition}

\begin{lemma}
Let $C$ be a congruence on the sesquiad $A$.
Among all additions on the monoid $A/C$ making the projection $A\to A/C$ a morphism of sesquiads, there is a minimal one $\Sigma_{\min}$ with the property that for every addition $\Sigma$ on $A/C$ making $A\to A/C$ a morphism of sesquiads, there exists a unique morphism of sesquiads $(A/C,\Sigma_{\min})\to (A/C,\Sigma)$ making the diagram
$$
\xymatrix{
A \ar[r] \ar[dr] & \(A/C,\Sigma_{\min}\) \ar@{.>}[d]^{\exists !}\\
& (A/C,\Sigma)
}
$$
commutative.
This means that the dotted morphism is the identity map, so this is equivalent to saying that the sum $\Sigma$ is defined on a possibly larger set than $\Sigma_{\min}$.
\end{lemma}

\begin{proof}
Let $(\Sigma_i)_{i\in I}$ be the family of all additions on $A/C$ making the projection a morphism of sesquiads.
For $i\in I$ let $R_i=R_{(A/C,\Sigma_i)}$ be the corresponding universal ring.
Let $J\subset R_A$ be the kernel of the map
$$
R_A\ \to\ \prod_{i\in I}R_i
$$
Then $A/C$ is a submonoid of $R_A/J$ and the addition on $A/C$ given by this inclusion $A/C\hookrightarrow R_A/J$ has the desired property.
\end{proof}

When we take the quotient $A/C$ of a sesquiad $A$ by a congruence $C$, we always install the minimal addition on $A/C$.

\begin{lemma}
\begin{enumerate}[\rm (a)]
\item For a congruence $C$ on a sesquiad $A$, let $J(C)$ denote the kernel of the map $R_A\to R_{A/C}$.
Then $J(C)$ is the ideal in $R_A$ generated by all elements of the form $a-b$, where $a\sim_C b$.
\item Let $I(A)$ be the kernel of $\Z_0A\to R_A$. Then the induced map $I(A)\to I(A/C)$ is surjective.
\end{enumerate}
\end{lemma}

\begin{proof}
(a) Let $\a\subset R_A$ be the ideal generated by all $a-b$ with $a\sim_Cb$.
The monoid $A/C$ maps to $R_A/\a$, and since $\a\subset J(C)$, the map $A/C\to R_A/\a$ is injective.
Hence $R_A/\a$ defines an addition on $A/C$ making the projection a morphism, therefore $R_{A/C}$ maps to $R_A/\a$ and the two are isomorphic.

(b) Consider the diagram with exact rows and columns,
$$
\xymatrix{
&&0\ar[d]&0\ar[d]\\
&&I(A)\ar[r]^\beta\ar[d]&I(A/C)\ar[d]\\
0\ar[r]&J_0(C)\ar[r]\ar[d]_\al&\Z_0A\ar[r]\ar[d]&\Z_0A/C\ar[r]\ar[d]&0\\
0\ar[r]&K\ar[r]&R_A\ar[r]\ar[d]&R_{A/C}\ar[r]\ar[d]&0\\
&&0&0
}
$$
Part (a) says that $J_0(C)$ and $K$ are both generated by all $a-b$ with $a\sim_Cb$, hence the map $\al$ is surjective.
From here it is a matter of diagram chase to show that $\beta$ is surjective, as well.
\end{proof}

If $\ph:A\to B$ is a sesquiad morphism, then we define its \e{congruence kernel} to be the congruence
$$
E=\ker_c(\ph)
$$
given by
$$
a\sim_E a'\qLR \ph(a)=\ph(a').
$$
This is indeed a congruence, as a compatible additive structure on $A/E$ is given by $A/E\hookrightarrow R_B$.
Note, however, that the addition induced on $A/E$ by the inclusion in $R_B$ is in general richer than the minimal addition of $A/E$. 

\begin{lemma}
Let $(C_k)$ be a family of congruences on $A$.
Then their intersection $C=\bigcap_kC_k$, as a subset of $A\times A$, is a congruence again.
\end{lemma}

\begin{proof} $C$ is the congruence kernel of the homomorphism
$
C\to\prod_k {A/C_k}.
$
\end{proof}

\begin{definition}
We say a sesquiad $A$ has \e{no zero divisors}, if
$$
ab=0\qR \(a=0\ \ \ \text{or}\ \ \ b=0\).
$$
We call a sesquiad \e{integral}, if $1\ne 0$ and
$$
af=bf\qR \(a=b\ \ \ \text{or}\ \ \ f=0\).
$$
An integral sesquiad has no zero divisors, but the converse does not hold in general.
A subsesquiad of an integral domain is an integral sesquiad.
\end{definition}

If $A$ is a sesquiad and $I$ is an ideal of the ring $R_A$, we also write $A/I$ for the sesquiad, which is the image of $A$ in $R_A/I$.
If $\a\subset A$ is an ideal of the sesquiad $A$, then we apply this to the ideal $I=(\a)$ of the ring $R_A$, so we write
$$
A/\a=A/(\a)=A/\ker_c(A\to R_{A}/(\a)).
$$
Note that if $A$ is a monoid, i.e., the addition is trivial, then $A/\a$ is the monoid, which is obtained by collapsing $\a$ to one point.

\begin{corollary}
Let $\a$ be an ideal in the sesquiad $A$.
Then the universal ring of $A/\a$ is $R_A/(\a)$.
\end{corollary}

\begin{proof}
The commutativity of the diagram
$$
\xymatrix{
A\ar[r]\ar[d]&R_A\ar[d]\\ A/\a\ar[r]&R_{A/\a}
}
$$
implies that $R_{A/\a}$ is a quotient of $R_A$.
The ideal $(\a)$ is the smallest such that $A/\a$ maps to the quotient, therefore the corollary.
\end{proof}

\subsection{Tensor product}
We define the tensor product $A\otimes B$ of two sesquiads $A$ and $B$ to be the submonoid of $R_A\otimes_\Z R_B$ consisting of all elements of the form $a\otimes b$ with $a\in A$ and $b\in B$.
We equip it with the additive structure of the ring $R_A\otimes R_B$.

\begin{lemma}
\begin{enumerate}[\rm (a)]
\item For any three sesquiads $A,B,C$ there is a functorial isomorphism
$$
\Hom(A,C)\times\Hom(B,C)\ \cong\ \Hom(A\otimes B,C).
$$
\item The universal ring of $A\otimes B$ is $R_A\otimes R_B$.
\end{enumerate}
\end{lemma}

\begin{proof}
(a) Let $\al:A\to C$ and $\be:B\to C$ be sesquiad morphisms. Define $\al\otimes\be:A\otimes B\to C$ by
$$
\al\otimes\be(a\otimes b)=\al(a)\be(b).
$$
The map $(\al,\be)\mapsto \al\otimes\be$ is easily seen to be a functorial bijection.

(b) Let $\phi:A\otimes B\to R$ be a morphism to a ring $R$.
By (a) we get ring homomorphisms $R_A\to R$ and $R_B\to R$ and their tensor product gives $R_A\otimes R_B\to R$ factoring $\phi$.
As $A\otimes B$ generates the ring $R_A\otimes R_B$, the ring homomorphism $R_A\otimes R_B\to r$ is uniquely determined by $\phi$, whence the claim.
\end{proof}

\subsection{Alternative definition}\label{alternative}
In this section we give an alternative definition of a sesquiad.

\begin{definition}
A \e{monoidal pair} we mean a pair $(A,R)$ consisting of a ring $R$ and a multiplicative submonoid $A$ of $R$ such that $A$ contains zero and $R$ is generated as a ring by $A$.

A \e{morphism} of monoidal pairs $\ph:(A,R)\to(B,S)$ is a ring homomorphism $\ph:R\to S$ mapping $A$ to $B$, i.e., $\ph(A)\subset B$.
\end{definition}

\begin{proposition}
Mapping a sesquiad $(A,\Sigma)$ to the pair $(A,R_A)$ is an equivalence of categories from the category of sesquiads to the category of monoidal pairs.
The inverse is given by mapping a pair $(A,R)$ to $(A,\Sigma)$ where $\Sigma$ is the addition given by the inclusion $A\subset R$.
\end{proposition}

\begin{proof}
The only nontrivial part is this: given a pair $(A,R)$, let $R_A$ be the universal ring of the sesquiad structure given by $A\subset R$. Then the natural map $R_A\to R$ is an isomorphism of rings, as follows from Proposition \ref{prop1.2.7} together with the fact that $R$ is generated by $A$.
\end{proof}

So a monoidal pair might be considered an equivalent definition of a sesquiad.
Likewise, one defnes a \e{semi-monoidal pair} $(A,S)$ consisting of a semi-ring $S$ together with a monoid $A\subset S$ which generates $S$ as a semi-ring.
These objects are called \e{blueprints} in \cite{blueprints}.
One could also  call them \e{semi-sesquiads}.

%: 
\section{Spectrum}
\subsection{Definition}
Let $A$ be  a sesquiad.
A congruence $\CC$ is called \e{prime} if $A/\CC$  is integral.
This is equivalent to $1\nsim_\CC 0$ and
$$
(af,bf)\in\CC\quad\Leftrightarrow\quad (a,b)\in\CC\text{  or  } (f,0)\in\CC.
$$

\begin{lemma}
Every congruence $\CC\ne A\times A$ is contained in a prime congruence.
\end{lemma}

\begin{proof}
This is equivalent to showing that $B=A/\CC$ has a prime congruence.
Let $\ph:B\to R$ be a representing morphism and let $\p$ be a prime ideal in the ring $R$.
Then $E=\ker_c(B\to R/\p)$ is a prime congruence on $B$.\end{proof}

\begin{definition}
Let $\spec_c A$ denote the set of all prime congruences with the topology generated by all sets of the form
$$
D(a,b)=\{\CC\in\spec_c A: (a,b)\notin\CC\}, \qquad a,b\in A.
$$
\end{definition}

\begin{theorem}
The space $\spec_c A$ is compact.
\end{theorem}

\begin{proof}
As the topology of $\spec_c A$ is generated by the sets $D(a,b)$,
the Alexander subbase theorem implies that it suffices to show that a given cover $\spec_cA=\bigcup_{i\in I}D(a_i,b_i)$ admits a finite subcover.
Taking complements we define
$$
C(a,b)=D(a,b)^c=\{E\in\spec_cA: a\sim_E b\}.
$$
We assume that for every finite set $F\subset I$ the intersection $\bigcap_{i\in F}C(a_i,b_i)$ is non-empty and we have to show that $\bigcap_{i\in I}C(a_i,b_i)\ne \emptyset$.
For $i\in I$ let $\a_i$ be the ideal in $R_A$ generated by $a_i-b_i$.
Let $\a$ be the ideal generated by all $\a_i$, so $\a=\sum_i\a_i$.
We claim that $\a$ is not the whole ring, i.e., that $1\notin\a$.
For assume that $1\in\a$, then there exists a finite set $F\subset I$ with $1\in\sum_{i\in F}\a_i$.
By our assumption, there exists a prime congruence $E\in \bigcap_{i\in F}C(a_i,b_i)$.
The ring homomorphism $P:R_A\to R_{A/E}$ annihilates the ideal $\sum_{i\in F}\a_i$, which is the whole ring, a contradiction!
It follows that the ring $R_A/\a$ is not the zero-ring, hence possesses a prime ideal $\p$.
The congruence 
$
\ker_c\(A\to (R_A/\a)/\p\)
$
then lies in $\bigcap_{i\in I}C(a_i,b_i)$, which therefore is non-empty.
\end{proof}

An ideal $J$ of the ring $R_A$ defines a congruence $C=C(J)$ by $a\sim_Cb\ \Leftrightarrow\ a-b\in J$.
In the converse direction, a congruence $C$ defines an ideal $I(C)$ of the ring $R_A$, which is the kernel of the ring homomorphism $R_A\to R_{A/C}$. 

\begin{lemma}\label{lem2.1.4}
The maps $C:\{ \text{ideals of }R_A\}\to\{\text{congruences on }A\}$ and $I$ in the opposite direction have the following properties:
\begin{enumerate}[\rm (a)]
\item $C\circ I=\Id$, so in particular, $C$ is surjective and $I$ is injective.
\item $C$ is  compatible with intersections,
$$
C\(\bigcap_kJ_k\)=\bigcap_k C(J_k).
$$
\item $C$ and $I$ are monotonic, i.e., if $J_1\subset J_2$, then $C(J_1)\subset C(J_2)$ and likewise for $I$.
\item If $J$ is a prime ideal, then $C(J)$ is a prime congruence, so $C$ induces a map $\spec R_A\to\spec_cA$. This map is continuous.
\item The ideal $I(E)$ is generated by all elements of the form $(a-b)$ where $a\sim_E b$.
\end{enumerate}
\end{lemma}

\begin{proof}
The proofs of (a) to (d) are straightforward.
For (e) one considers the ideal $J$ generated by all $(a-b)$ with $a\sim_Eb$. Then one sees that the ring $R_A/J$ is the universal ring of $A/E$.
\end{proof}

\begin{definition}
A sesquiad morphism $\ph:A\to B$ induces a continuous map
$$
\ph^*:\spec_c B\to\spec_cA
$$
given by
$$
E\mapsto \ker_c(A\to B/E).
$$

Recall the \e{max-spectrum} of a ring $R$, which is the set of all maximal ideals of $R$ and is denoted by $\Mspec R$.
\end{definition}

\begin{lemma}\label{meetsevery}
Let $A$ be a sesquiad and let $\rho:A\to R_A$ be the canonical morphism.
The image of $\rho^*:\spec R_A\to\spec_cA$ meets every non-empty closed set in $\spec_cA$.
Even its restriction to the max-spectrum $\Mspec R_A\to\spec_cA$ has the same property.
\end{lemma}

\begin{proof}
Let $E\in\spec_cA$.
We show that $\rho^*(\Mspec R_A)$ meets the closure $\ol{\{E\}}$ of $E$.
The kernel $I$ of $R_A\to R_{A/E}$ is a proper ideal.
It therefore lies in some maximal ideal $m$ of $R_A$.
The congruence $\ker_c(A\to R_A/m)$ lies in the closure  of $E$ and in $\rho^*(\Mspec R_A)$.
\end{proof}

\begin{lemma}
Let $A,B$ be sesquiads.
The two projections of the sesquiad $A\times B$ induce a homeomorphism
$$
\phi:\spec_cA\cupdot\spec_cB\tto\cong\spec_c(A\times B).
$$
\end{lemma}

\begin{proof}
For $E\in\spec_cA$ we define $\phi(E)$ as the congruence kernel of $A\times B\to A/E$ and likewise for $B$.
It is clear that the map $\phi$ is injective.
For surjectivity note that because of $(1,0)(0,1)=0$ every $E\in\spec_c(A\times B)$ has to send one of these two elements to zero.
By $(a,b)=(a,0)+(0,b)$, the class of $(a,b)$ depends only on one of the entries. This ensures surjectivity.
The continuity of $\phi$ and its inverse are easy.
\end{proof}

\begin{definition}
We also consider the \e{Zariski spectrum} $\spec_zA$, which is the set of all prime ideals of the sesquiad $A$.
It is equipped with the \e{Zariski topology}, which is generated by all sets of the form
$$
D(f)=\{\p\in\spec_zA: f\notin\p\},
$$
for $f\in A$.
\end{definition}

There is a canonical map
$$
\CZ:\spec_cA\to\spec_zA
$$
sending a prime congruence $x$ to the class $[0]_x$ of the zero element.
This map is continuous.
It is called the \e{zero class map}.

\begin{lemma}
\begin{enumerate}[\rm (a)]
\item Let $R$ be a ring and $S\subset R\sm\{0\}$ be a multiplicative monoid.
Then there exists a prime ideal $\p$ in $R$ such that $\p\cap S=\emptyset$.
\item The zero class map is surjective.
\end{enumerate}
\end{lemma}

\begin{proof}
(a) By the Lemma of Zorn we get an ideal $\p$ of $R$ which is maximal with the property that $\p\cap S=\emptyset$.
We show that it is prime.
Let $x,y\in R$ with $xy\in\p$ and $x\notin\p$.
Then the ideal $\sp{\p,x}$ generated by $\p$ and $x$ has non-empty intersection with $S$, so there is $s\in S$ with $s=p+rx$ for some $p\in\p$ and $r\in R$.
Then $sy=py+rxy$ lies in $\p$.
Now \emph{assume} $y\notin \p$, then $\sp{y,\p}$ has non-empty intersection with $S$, so there is $t\in S$ with $t=p'+yr'$ for some $p'\in\p$ and $r'\in R$.
Then $st=sp'+syr$ lies in $\p$ and in $S$, a contradiction!
Hence it follows $y\in \p$, so $\p$ is prime.

(b) Let $\a$ be a prime ideal of the sesquiad $A$ and let $P$ be the ideal of $R_A$ generated by $\p$.
Let $S$ be the image of $A$ in the ring $R-R_A/P$.
Then $S\subset R\sm\{0\}$ is a monoid.
By part (a) there exists a prime ideal $\p$ of $R$ with $\p\cap S=\emptyset$.
Let $C$ be the congruence $\ker(A\to R/\p)$,
then $\CZ(C)=\a$.
\end{proof}

\subsection{Nilradical}
The \e{nilradical} of the sesquiad $A$ is defined to be the intersection of all prime congruences, so
$$
\Nil(A)=\bigcap_{E\in\spec_cA}E=\ker_c\(A\to\prod_{E\in\spec_cA}A/E\).
$$

\begin{proposition}
If a pair $(a,b)$ lies in the nilradical, then $(a-b)$ is nilpotent in the ring $R_A$.
So we have
$$
(a,b)\in\Nil(A)\quad\Rightarrow\quad (a-b)\in\Nil(R_A),
$$
where for a ring $R$ by $\Nil(R)$ we denote its set of all nilpotent elements.
\end{proposition}

\begin{proof}
Let $E\in\spec_cA$, let $\q$ be a prime ideal of the ring $R_{A/E}$ and let $\p=\ker(R_A\to R_{A/E}/\q)$.
Let $E_\p$ be the congruence kernel of $A\to R_A/\p$, then $E_\p\supset E$ and therefore
$$
\bigcap_{E\in\spec_cA}E\subset\bigcap_{\p\in\spec R_A}E_\p=\ker_c\(A\to\prod_{\p\in\spec R_A}R_A/\p\).
$$
It follows that $(a,b)\in\Nil(A)$ implies  $(a-b)\in\bigcap_{\p\in\spec R}\p=\Nil(R_A)$.
\phantom{ }\end{proof}

\begin{example}
The following example shows that in the preceeding proposition the converse direction is false in general.
Consider the sesquiad $A=\{0,1,3\}\subset\Z/4\Z$ with the addition of the ring $\Z/4\Z$.
Then with $a=3$ and $b=1$ the element $a-b=2$ is nilpotent in $R_A=\Z/4\Z$, but, as $A$ is integral, we have $\Delta\in\spec_cA$ and so $(a,b)\notin\Nil(A)$.
\end{example}

\begin{definition}
We say that a sesquiad $A$ is \e{reduced}, if $\Nil(A)=\Delta$. 
\end{definition}

\begin{lemma}\label{lem2.2.4}
For any sesquiad $A$, the quotient $A/\Nil(A)$ is reduced. It is the largest reduced quotient of $A$, we call it the \e{reduction} of $A$ and write
$$
 A^\red=A/\Nil(A).
 $$
The map $A\to A^\red$ induces a homeomorphism 
$$
\spec_c A^\red\tto\cong\spec_cA.
$$
\end{lemma}

\begin{proof}
Let $a,b\in A$ with $\bar a, \bar b$ their classes in $A/\Nil(A)$ and assume $(\bar a,\bar b)\in\Nil(A^\red)$.
Any $E\in\spec_cA$ contains $\Nil(A)$ and therefore induces an element of $\spec_cA^\red$.
Therefore $(a,b)\in\Nil(A)$, whence the first claim.
The rest is clear.
\end{proof}

\subsection{Irreducible and noetherian spaces}
Recall that a topological space $X$ is called \e{irreducible} if $X=A\cup B$ for some closed sets $A,B$ implies that $A=X$ or $B=X$.
This is equivalent to saying that any two non-empty open sets have a non-empty intersection.
This again is equivalent to 
saying that any non-empty open subset is dense.

\begin{definition}
A topological space $X$ is called \e{sober}, if every non-empty irreducible closed subset has a unique generic point.
\end{definition}

If $I$ is an ideal of the ring $R_A$, we write $A/I$ for the sesquiad which is the image of $A$ in $R_A/I$.

\begin{proposition}
For a sesquiad $A$ we have
\begin{enumerate}[\rm (a)]
\item the space $X=\spec_cA$ is irreducible if and only if $A^\red$ is integral
\item the space $\spec_cA$ is sober.
\end{enumerate}
\end{proposition}

\begin{proof}
For (a) suppose that $\spec_cA$ is irreducible.
By Lemma \ref{lem2.2.4} we can replace $A$ with $A^\red$ and assume that $A$ is reduced.
We show that $A$ is integral.
Let $fx=fy$ hold in $A$ with $f\ne 0$.
This implies that $D(f,0)\cap D(x,y)=\emptyset$.
As $A$ is reduced, $f\ne 0$ implies $D(f,0)\ne \emptyset$, therefore, as $X$ is irreducible, we have $D(x,y)=\emptyset$, therefore, again as $A$ is reduced, $x=y$, which means that $A$ is integral.

For the converse direction assume $A$ is integral.
Then $\Delta$ lies in every set of the form $D(a,b)$ with $a\ne b$, so $\Delta$ lies in every open set.
So any two non-empty open sets have a non-empty intersection.

(b) Let $S\subset \spec_cA$ be an irreducible, non-empty, closed subset.
Consider the congruence
$$
\eta=\bigcap_{E\in S}E.
$$
We claim that $\eta$ is the desired generic point.
Since $\eta$ is contained in every $E\in S$, which is equivalent to $E\in\ol\eta$, the proof of the claim reduces to showing that $\eta\in S$.

To start with, we need to know that $\eta$ is prime.
So let $f,x,y\in A$ with $fx\sim_\eta fy$.
Then the same equivalence holds under all $E\in S$.
Since every $E\in S$ is prime, one has $f\sim_E 0$ or $x\sim_E y$.
This means that $S$ is the union of the two closed sets $C(f,0)\cup S$ and $C(x,y)\cup S$.
Therefore, one of them equals $S$.
If $S\subset C(f,0)$, then $f\sim_\eta 0$ and if $S\subset C(x,y)$, then $x\sim_\eta y$, so $\eta$ is prime indeed.
To finally show that $\eta\in S$ we show that $S$ has non-empty intersection with every neighborhood of $\eta$. As $S$ is closed, this implies $\eta\in S$.
A neighborhood base for $\eta$ is given by the sets of the form $D(a,b)$, $a,b\in A^n$ which contain $\eta$.
So assume $\eta\in D(a,b)=D(a_1,b_1)\cap\dots\cap D(a_n,b_n)$.
This implies that for each $j$ there is $E\in S$ with $a_j\nsim_E b_j$, so $D(a_j,b_j)\cap S\ne\emptyset$.
As $S$ is irreducible, it follows that $D(a,b)\cap S\ne \emptyset$ as promised.
\end{proof}

Recall that a topological space $X$ is said to be \e{noetherian}, if every descending sequence of closed sets is eventually stationary.
Equivalently, every ascending sequence of open sets is eventually stationary.
Closed or open subsets of a noetherian space are noetherian.
A space is noetherian if and only if every open subset is compact.

\begin{definition}
A sesquiad $A$ is called \e{noetherian}, if every congruence $C$ on $A$ is finitely generated, which means that there are $a,b\in A^n$ such that $C$ is generated by $a_1\sim b_1,\dots,a_n\sim b_n$. 
This again means that $C$ is the intersection of all congruences $E$ with $a_1\sim_E b_1,\dots,a_n\sim_E b_n$.
\end{definition}

\begin{lemma}
Let $A$ be  a sesquiad.
\begin{enumerate}[\rm (a)]
\item $A$ is noetherian if and only if every ascending sequence $C_1\subset C_2\subset\dots$ of congruences is eventually stationary.
\item If $A$ is noetherian, and $\phi:A\twoheadrightarrow B$ is a surjective morphism of sesquiads, then $B$ is noetherian.
\item Every finitely generated sesquiad is noetherian.
\item If $A$ is noetherian, then the topological space $\spec_cA$ is noetherian.
\end{enumerate}
\end{lemma}

\begin{proof}
(a) Let $A$ be noetherian and $C_1\subset C_2\subset\dots$ be an ascending sequence of congruences.
Then $C=\bigcup_{j=1}^\infty C_j$ is a congruence.
This is not completely trivial, as a representing ring has to be found.
This is
$$
R=\lim_{\stack{\to}j} R_{A/C_j}.
$$
As $A$ is noetherian, the congruence $C$ is finitely generated. So the generating relations are present at some finite stage.
This proves the ascending sequence condition.
For the converse direction assume the ascending sequence condition and let $C$ be a congruence on $A$.
Assume that $C$ is not finitely generated.
Then there exists a sequence $(a_j,b_j)\in A\times A$ auch that the congruences $C_n$ generated by $a_1\sim b_1,\dots,a_n\sim b_n$ are all distinct.
This contradicts the ascending sequence condition.

(b) Let $E$ be the congruence kernel of $\phi$.
Then $B$ coincides with $A/E$ as a monoid, but the addition on $B$ may be richer.
In any case, congruences on $B$ are given by congruences on $A$ which contain $E$.
Now use the ascending sequence condition of (a).

(c)
Using (b), it suffices to consider the sesquiad $A=\sp{T_1,\dots,T_n}$ freely generated by $T_1,\dots,T_n$.
The universal ring is the polynomial ring$R_A=\Z[T_1,\dots,T_n]$, which is noetherian.
Let $E$ be a congruence on $A$ and let $I$ be the ideal of $R_A$ generated by all $a-b$ with $a\sim_E b$.
Since the ring $\Z[T_1,\dots,T_n]$ is noetherian, the ideal $I$ is generated by finitely many 
$a_1-b_1,\dots,a_n-b_n$.
As the congruence $E$ is given by $I$, the congruence $E$ is generated by $a_1\sim b_1,\dots,a_n\sim b_n$.

(d) 
Let $C_1\supset C_2\supset\dots$ be a descending sequence of closed sets.
As every closed set is an intersection of sets of the form $C(a,b)=C(a_1,b_1)\cup\dots\cup C(a_n,b_n)$ for some $a,b\in A^n$ it suffices to prove the descending sequence condition for closed  sets of the form
$C_i=C(a^i,b^i)$ with $a^i,b^i\in A^{n_i}$.
Then it is enough to prove it for closed sets of the form $C_i=C(a_i,b_i)$ with $a_i,b_i\in A$.
In this situation, let $c_i$ be the congruence generated by $a_1\sim b_1,\dots,a_i\sim b_i$.
The sequence $c_i$ is eventually stationary, so then is $C_i$.
\end{proof}

\subsection{Structure sheaf}

\begin{lemma}
Let $\ph:A\to B$ be a surjective morphism of sesquiads, then the map $\ph^*:\spec_c B\to\spec_cA$ is a homeomorphism onto its image, which is a closed subset of $\spec_cA$.
\end{lemma}

\begin{proof}
It is a continuous map which clearly is injective.
Let $c$ be the congruence $\ker_c(\ph)$.
Then the image of $\ph$ is the intersection of the following closed subsets.
Firstly the closure $\ol c$, i.e., the set of all $E\in\spec_cA$ with $E\supset c$.
Next for all $a\in A^n$, $k\in\Z^n$ and $\al\in A$ such that $\sum_{j=1}^nk_j\ph(a_j)=\ph(\al)$, the set of all $E\in\spec_cA$ that contain
$$
\ker_c\(A\to R_A\left/\(\sum_{j=1}^nk_ja_j-\al\)\right.\).
$$ 
It follows that the image is closed.
Finally, $\ph^*$ is a homeomorphism onto its image, since
\begin{align*}
\ph^*(D_B(b,b'))&=\left[\bigcup_{\stack{\ph(a)=b}{\ph(a')=b'}}D_A(a,a')\right]\cap \Im(\ph^*).
\tag*{$\square$}
\end{align*}

For a prime congruence $E$ on a sesquiad $A$ we define
$$
S_E=\sp{(a-b):0\nsim_Ea\nsim_Eb\nsim_E1}
$$
as the submonoid of $R_A$ generated by all elements of the form $(a-b)$ where $a,b\in A$ with $0\nsim_E a\nsim_E b\nsim 1$.
We define the localization
$$
A_E=S_E^{-1}A
$$
as the subsesquiad of $S_E^{-1}R_A$ generated by $A$ and $S_E^{-1}$ with the addition induced by the inclusion into the localized ring $S_E^{-1}R_A$.

Note that by the equation
$$
1+\frac b{a-b}=\frac a{a-b},
$$
the sesquiad $A_E$ can have nontrivial addition even if the addition on $A$ is trivial.

\begin{lemma}
\begin{enumerate}[\rm (a)]
\item The universal ring of $A_E$ is $S_E^{-1} R_A$.
\item For every  morphism of sesquiads $\phi:A\to B$ with $\phi_R(S_E)\subset B^\times$ there exists a unique morphism of sesquiads making the diagram
$$
\xymatrix{
A \ar[r]^\phi \ar[rd] & B\\
& A_E\ar@{.>}[u]
}
$$
commutative.
Here $\phi_R$ is the induced ring homomorphism $R_A\to R_B$.
\end{enumerate}\end{lemma}

\begin{proof}
(a)
Every sesquiad morphism $S_E^{-1}A\to R$ to a ring $R$  induces by pullback a sesquiad morphism $A\to R$ which factors uniquely over a ring morphism $\ph:R_A\to R$.
Then $\ph$ maps $S_E$ into the unit group and so factors uniquely over $S_E^{-1}R_A$.
The diagram
$$
\xymatrix{
S_E^{-1}A\ar[r]\ar[rd] & R\\
& S_E^{-1}R_A\ar@{.>}[u]
}
$$
commutes, so part (a) follows.
For (b) let $\phi$ and $\phi_R$ be as above.
By the universal property of the ring localization, $\phi_R$ factors uniquely over $S_E^{-1}R_A$. This gives the claim.
\end{proof}

\begin{lemma}
Let $A$ be a sesquiad, then the common congruence kernel of all localizations is the trivial congruence, i.e.,
$$
\bigcap_{p\in\spec_cA}\ker_c(A\to A_p)=\Delta.
$$
\end{lemma}

\begin{proof}
We show more sharply that
$$
\bigcap_{p\in\spec R_A}\ker_c(A\to R_{A,p})=\Delta.
$$
This follows from the general fact that for any ring $R$ one has
$$
\bigcap_{p\in\spec R}\ker(R\to R_p)=\{ 0\}.
$$
To prove this, let $x\in R$ be non-zero. 
Then there exists a prime ideal $p$ containing the annihilator $\{ y\in R: xy=0\}$ of $x$.
It follows that $x$ is nonzero in $R_p$.
\end{proof}

\begin{definition}
For an open set $U\subset\spec_c A$ a \e{section} is a map
$$
s:U\to\coprod_{E\in U}A_E
$$
with $s(E)\in A_E$ for every $E\in U$, such that $s$ is locally a quotient of elements of $R_A$, i.e., for every $E\in U$ there exists an open set $V$ with $E\in V\subset U$ and
$$
a\in A,\qquad f\in\bigcap_{F\in V}S_F,
$$
such that for every $F\in V$ one has
$$
s(F)=\frac af.
$$
We write $\CO(U)$ for the set of sections over $U$.
It is clear that this set forms a monoid.
We claim that it inherits a natural structure of a sesquiad.
Recall that $R_{A_E}=S_E^{-1}R_A$ and 
let $\CR(U)$ denote the set of all maps
$$
s:U\to\coprod_{E\in U}S_E^{-1}R_A
$$
which are locally, on $V\subset U$ say, of the form $s(F)=\frac \al f$ with 
$$
\al\in R_A,\qquad f\in\bigcap_{F\in V}S_F.
$$
Then $\CR$ is a ring-valued sheaf on $\spec_cA$ and the inclusion of  the monoid $\CO(U)\hookrightarrow \CR(U)$ makes $\CO(U)$ a sesquiad.
\end{definition}

\begin{theorem}\label{thm2.3.8}
Let $A$ be a sesquiad.
\begin{enumerate}[\rm (a)]
\item The sections form a sheaf $\CO$ of sesquiads on $\spec_c A$.
\item Point evaluation gives an isomorphism $\CO_E\tto\cong A_E$, where $\CO_E$ is the stalk of $\CO$ at $E\in\spec_c A$.
\item The sesquiad $\Ga A=\Ga(\CO_A)$ of global sections contains $A$ as a subsesquiad and is contained in $R_A$, i.e.,
$$
A\hookrightarrow\Ga A\hookrightarrow R_A.
$$
The universal ring of $\Ga A$ is the ring $R_A$.
\item For every morphism $\ph:A\to B$ of sesquiads there is a unique morphism $\ph_\Ga:\Ga A\to\Ga B$ making the diagram
$$
\xymatrix{
A\ar[r]^\ph\ar[d] & B\ar[d]\\
\Ga A\ar[r]^{\ph_\Ga}& \Ga B
}
$$
commutative.
\item For sesquiads $A$ and $B$ let $\Hom_{A,B}(\Ga A,\Ga B)$ denote the set of all morphisms of sesquiads $\al:\Ga A\to \Ga B$ with $\al(A)\subset B$.\\
The map $\ph\mapsto\ph_\Ga$ of the previous part is an ismorphism of bifunctors
$$
\psi:\Hom(A,B)\tto\cong\Hom_{A,B}(\Ga A,\Ga B).
$$
\end{enumerate}
\end{theorem}

\begin{proof}
The parts (a) and (b) are clear.
We start the proof of (c) and (d) with the construction of $\ph_\Ga$ in part (d).
Let $s:\spec_cA\to \coprod_{E\in\spec_cA}A_E$ be a global section.
Then $s\circ\ph^*$ is a map from $\spec_cB$ to 
$\coprod_{F\in\spec_cB}A_{\ph^*F}$.
For every $F\in\spec_cB$ the morphism $\ph$ localizes to a morphism $\ph_F:A_{\ph^*F}\to B_F$.
We apply this localization after $s\circ\ph^*$ to obtain
$$
\ph_\Ga(s):\spec_cB\to\coprod_{F\in\spec_cB}B_F.
$$
A local representation $s(E)=\frac af$ is mapped to a local representation $\frac{\ph(a)}{\ph(f)}$, so $\ph_\Ga(s)$ is indeed a section.
Elements of $A$ are mapped to their images under $\ph$, so $\ph_\Ga$ makes the diagram commutative.

We apply this construction in the special case of $B$ being $R_A$ to obtain the commutative diagram
$$
\xymatrix{
A\ar[r]^\rho\ar[d]^\ga & R_A\ar[d]^\cong\\
\Ga A\ar[r]^{\rho_\Ga} & \Ga(\CO_{R_A}).
}
$$
We claim that $\rho_\Ga$ is injective.
First consider the set $X=\rho^*(\spec(R_A))$.
By Lemma \ref{meetsevery} we know that $X$ meets every closed set in $\spec_cA$, therefore, any open set that contains $X$, must be all of $\spec_cA$.

Let $s,t$ be two global sections of $\CO_A$ with $\rho_\Ga(s)=\rho_\Ga(t)$.
We show that the two sections coincide on the image $X=\rho^*(\spec(R_A))$ in $\spec_cA$.
For this let $\p\in\spec(R_A)$.
Note that the localization $\rho_\p:A_{\rho^*\p}\to R_{A,\p}$ factorizes over the universal ring $R_{A_{\rho^*\p}}=S_{\rho^*\p}^{-1}R_A$.
Let $\q$ be a prime ideal of the ring $S_{\rho^*\p}^{-1}R_A$ and let $\q'$ be its pullback to $R_A$.
Then $S_{q'}\supset S_{\rho^*\p}$ and the localization map induces an isomorphism
$$
\(S_{\rho^*\p}^{-1} R_A\)_q\cong R_{A,q'}.
$$
Since $\rho_\Ga(s)=\rho_\Ga(t)$, these two sections coincide at every $\q'$ and therefore they coincide in $S_{\rho^*\p}^{-1} R_A$, hence they coincide on $X$.
The set $\{s=t\}\subset\spec_cA$ is open and contains 
$X$, so it must be all of $\spec_cA$, which means $s=t$.
This shows the claimed injectivity of $\rho_\Ga$, which in turn shows the injectivity of $\ga$ as well.
By Lemma \ref{1.2.3}, $\rho_\Ga$ induces an identification $R_A\cong R_{\Ga A}$.
This in particular implies that $\rho_\Ga$ is the unique sesquiad morphism making the diagram
$$
\xymatrix{
A\ar[r]^\rho\ar[dr]_\ga & R_A\\
& \Ga A\ar[u]_{\rho_\Ga}
}
$$
commutative.

We now show the uniqueness of $\ph_\Ga$ in part (d).
Let $\al:\Ga A\to\Ga B$ be a second morphism making the diagram commutative.
We extend the diagram to
$$
\xymatrix{
A\ar[r]^\ph\ar[d] & B\ar[d]\\
\Ga A\ar[r]^\al \ar[d] &\Ga B\ar[d]\\
R_A=R_{\Ga A}\ar[r] & R_B=R_{\Ga B}
}
$$
Both squares commute. The commutativity of the whole diagram forces the lower horizontal arrow to be the unique one induced by $\ph$.
This implies $\al=\ph_\Ga$.

(e) Let $R:\Hom_{A,B}(\Ga A,\Ga B)$ be the restriction homomorphism $R(\al)=\al|_A$.
Then both $\psi$ and $R$ are functorial and $R$ is a left-inverse to $\psi$, so $\psi$ is injective and $R$ surjective.
All that we need to show is that $R$ is injective.
So let $\al,\be\in\Hom_{A,B}(\Ga A,\Ga B)$ with $R(\al)=R(\be)$, so $\al$ and $\be$  agree on $A$.
For every $F\in\spec_cB$ we thus have 
$F^*\df\al^*F=\be^*F$ and the localizations $\al_F,\be_F:A_{F^*}\to B_F$ agree.
The point-evaluation morphisms give a commutative diagram
$$
\xymatrix{
\Ga A\ar[r]^{\al,\be}\ar[d] &\Ga B\ar[d]\\
\displaystyle\prod_{F\in\spec_cB}A_{F^*}\ar[r] & \displaystyle\prod_{F\in\spec_cB}B_F
}
$$
which commutes for both, $\al$ or $\beta$ in the top row.
The injectivity of the right vertical arrow implies $\al=\be$.
\end{proof}

\begin{definition}
We call a sesquiad $A$ \e{conservative}, if the map $\ga: A\to\Ga A$ is an isomorphism  of sesquiads.
By the theorem, this is equivalent to $\ga$ being surjective.
\end{definition}

\begin{examples}\label{Ex2.3.10}
\item This is an example of an integral sesquiad which is not conservative:
Let $\F_7=\Z/7\Z$ be the field of seven elements.
The sesquiad $A=\{0,1,2,4\}\subset \F_7$ with the addition from $\F_7$ is not conservative, as $R_A=\F_7$ and $\spec_cA=\{\Delta\}$, so $\Ga A\cong A_\Delta$, which is all of $\F_7$.
\item We give an example of a non-conservative sesquiad with trivial addition:
Let $A=\{ 0,1,e\}$ with $e^2=e$ and trivial addition.
Then $\spec_cA=\{ e\sim 0,e\sim 1\}=\{E_0,E_1\}$.
One has $A_{E_0}=\{0,1\}$ and $A_{E_1}=\{0,1\}$.
As $\spec_cA$ has the discrete topology, we have $\Ga A=A_{E_0}\times A_{E_1}$, so $A$ is not conservative.
Moreover, $A$ is a monoid, i.e., has trivial addition, but $\Ga A$ has not, as the addition $(1,0)+(0,1)=(1,1)$ is defined in $\Ga A$
\item It may even happen that although $A$ is integral, the sesquiad $\Ga A$ is not.
An example for this is $A=\{0\}\cup(\Z/15\Z)^\times$, in which case $R_A=\Z/15\Z=\Z/3\Z\times\Z/ 5\Z$.
The spectrum has three points: $\Delta$, the congruence kernel of $A\to\Z/3\Z$ and the congruence kernel of $A\to\Z/5\Z$.
One gets $\Ga A=\Z/3\Z\times\Z/5\Z$ which is not integral.

In this example we also see that
$$
\spec_cA\ncong\spec_c\Ga A
$$
can happen, as the generic point of $\spec_cA$ is no longer present in $\spec_c\Ga A$.
\item We give an example of a sesquiad morphism $\ph:A\to B$ which is injective, but the induced morphism on sections $\ph_\Ga:\Ga A\to\Ga B$ is not injective.
Let $A=(\Z/6\Z)^\times\cup\{ 0\}=\{ 0,\pm 1\}$ with the addition from $\Z/6\Z$, and let $B=\Z/3\Z$.
Then the projection $A\to B$ has the claimed property.
\item The functor $\Ga$ also doesn't preserve surjectivity, as the following example shows.
Let $A=\{ 0,1,\tau,\tau^2,\dots\}$ be the free 
monoid with generator $\tau$.
Let $B=\{0,1,e\}$ be the monoid with the relation $e^2=e$ and let $\ph:A\to B$ be the monoid map with $\ph(\tau)=e$.
Then $\ph$ is surjective, but as $A$ is conservative and $B$ is not, the map $\ph_\Ga:\Ga A\to\Ga B$ is not surjective.
\end{examples}

\begin{proposition}\label{propcons}
\begin{enumerate}[\rm (a)]
\item \label{propcons1}A ring is a conservative sesquiad.
\item \label{propcons2}An integral sesquiad with trivial addition is conservative.
\item \label{propcons3}Let $A$ be a sesquiad and let $\ga:A\to\Ga A$ be the inclusion.
Suppose there is a morphism of sesquiads $\phi:\Ga A\to A$ such that $\phi\circ\ga=\Id_A$.
Then $A$ is conservative and $\phi$ is the identity.
\end{enumerate}
\end{proposition}

\begin{proof}
(\ref{propcons1}) is classical. For (\ref{propcons2}) let $A$ be integral with trivial addition.
Let $c$ be the congruence with two equivalence classes $[0]_c=\{0\}$ and $[1]_c= A_1= A\sm\{ 0\}$.
The monoid $A_1$ is cancellative and injects into its quotient group $G$.
The localization $A_c$ is the monoid $G_0=G\cup\{ 0\}$.
The stalk $\CR_c$ is the ring $S_c^{-1}R_A=S_c^{-1}\Z_0A=A^{-1}\Z_0A=\Z G$, which is the group ring and the point evaluation is an injective map from $R_A=\Z{A_1}$ into the stalk $\CR_c$.
So $\Ga A$, as a subset of $\Z G$, lies in the intersection of $A_c=G_0$ and $R_A$.
This intersection is $A$.

(\ref{propcons3})
Both maps $\ga$ and $\phi$ induce ring homomorphisms on $R_A$.
We get a commutative diagram
$$
\xymatrix{
R_A\ar[r]^{\ga_R\ \ \ }\ar[dr]_{(\phi\circ\ga)_R} & R_{\Ga A}=R_A\ar[d]^{\phi_R}\\
& R_A
}
$$
As $\phi\circ\ga=\Id$ and $\ga_R=\Id$ we get $\phi_R=\Id$. On the other hand, $\phi_R$ maps $\Ga A$ to $A$, therefore they are equal.
\end{proof}

\begin{example}
If $\ph:A\twoheadrightarrow B$ is a surjective morphism of sesquiads, then the map $\ph_\Ga:\Ga A\to\Ga B$ needs not be surjective as the following example shows.
Let $A=\{ 0,1,\tau,\tau^2,\dots\}$ be the free monoid generated by one element $\tau$, let $B=\{0,1,e\}$ with $e^2=e$ and let $\ph:A\to B$ sending $\tau$ to $e$.
Then $\Ga A=A$ and so the map $\ph_\Ga:\Ga A\to\Ga B$ factors over $B$.
As $B$ is not conservative, $\ph_\Ga$ is not surjective.
\end{example}

\begin{definition}
Let $n\in\N$.
For $a,b\in A^n$ we set
$$
D(a,b)=D(a_1,b_1)\cap\dots\cap D(a_n,b_n).
$$
Note that this includes the case $n=1$ and $(a,b)=(1,0)$, where $D(a,b)$ is the entire space $\spec_cA$.
\end{definition}

\begin{lemma}\label{lem2.3.10}
Let $A$ be a sesquiad.
The map $\rho:A\to R_A$ induces the map $\phi=\rho^*:\spec R_A\to\spec_cA$.
There is a natural isomorphism of sheaves 
$$
\phi_*\CO_{R_A}\tto\cong\CR.
$$
\end{lemma}

\begin{proof}
We compare the stalks of the two sheaves on the space $\spec_cA$.
For $E\in\spec_cA$ the stalk of $\CR$ is $\CR_E=R_{A_E}=S_E^{-1}R_A$.
The stalk of $\phi_*\CO_{R_A}$ is
\begin{align*}
\(\phi_*\CO_{R_A}\)_E&=\lim_{U\ni E}\CO_{R_A}(\phi^{-1}(U)).
\end{align*}
Since every open set $U$ is a union of sets of the form $D(a,b)$ with $a,b\in A^n$, it suffices to extend the limit over the latter.
We have $\CO_{R_A}(\phi^{-1}(D(a,b)))=\CO_{R_A}(D(f))=R_A[f^{-1}]=S_f^{-1}R_A$,
where $f=(a_1-b_1)\cdots(a_n-b_n)$.
We therefore get a map $\CO_{R_A}(\phi^{-1}(D(a,b)))\to\CR(D(a,b))$ which extends to a morphism of sheaves and we see that it is an isomorphism on stalks, therefore it is an isomorphism of sheaves, which proves the lemma.
\end{proof}

\subsection{$\Ga\Ga A=\Ga A$}
In this section we will prove the following theorem.

\begin{theorem}
Let $A$ be a sesquiad.
Then $\Ga A$ is conservative.
In other words, we have
$$
\Ga\Ga A=\Ga A.
$$
\end{theorem}

The proof will occupy the rest of the section.

\begin{definition}
A sesqiad $A\ne 0$ is called \e{simple}, if for every sesquiad $B\ne 0$, every sesquiad morphism $A\to B$ is injective.

A point $x\in\spec_cA$ is called a \e{closed point}, if $\{ x\}$ is a closed set.
\end{definition}

\begin{proposition}\label{prop2.5.3}
\begin{enumerate}[\rm (a)]
\item A simple sesquiad is integral.
\item For a sesquiad $F$ the following are equivalent:
\begin{enumerate}[\rm (i)]
\item $F$ is simple,
\item For every ideal $I\subset R_F$ with $I\ne R_F$, the sesquiad $F$ injects into $R_F/I$,
\item For every prime ideal $\p\subset R_F$ the sesquiad $F$ injects into $R_F/\p$,
\item For every maximal ideal $\m\subset R_F$ the sesquiad $F$ injects into $R_F/\m$.
\end{enumerate} 
\item A sesquiad $A$ is simple if and only if
$$
\spec_cA=\{\Delta\}.
$$
\item Let $X=\spec_cA$, where $A$ is a sesquiad.
An element $E\in X$ is a closed point if and only if the corresponding quotient $A/E$ is simple.
\item Every sesquiad $A\ne 0$ has a simple quotient $A/E$.
\item Every non-empty closed subset of $\spec_cA$ contains a closed point.
\item If $A$ is simple, then $\Ga A$ is integral.
\end{enumerate}
\end{proposition}

\begin{proof}
(a) We show that a non-integral sesquiad $A$ is not simple.
So let $f,x,y\in A$ with $f\ne 0$, $x\ne y$, and $fx=fy$.
Then $f$ is a zero-divisor in the ring $R_A$, so the ideal $fR_A$ is not the entire ring, i.e., map $A\to R_A/fR_A$ is non-trivial. 
As this map sends $f$ to zero, it is not injective as well. This means that $A$ is not simple.

(b) The only non-trivial assertion is (iv)$\Rightarrow$(i).
We decompose it in two steps. First we show (iv)$\Rightarrow$(iii) and then (iii)$\Rightarrow$(i).
Suppose (iv) holds and let $\p\subset R_F$ be a prime ideal.
Choose a maximal ideal $\m$ containing $\p$.
Consider the commutative diagram,
$$
\xymatrix{
F\ar[r]^\al\ar[rd]_\beta &R_F/\p\ar[d]\\
&R_F/\m.
}
$$
By assumption, $\beta$ is injective. By commutativity of the diagram, $\al$ must be injective, too.
Next assume (iii) holds and consider a sesquiad morphism $\ph:F\to B$.
This induces a ring homomorphism $\ph:R_F\to R_B$.
Let $J$ be a prime ideal of $R_B$, then $\ph^{-1}(J)$ is a prime ideal of $R_F$ so $F$ injects into $R_B/J$, hence $\ph$ is injective.

(c) Let $A$ be simple. As $A$ is integral, $\Delta$ lies in $\spec_cA$.
Let $C\in\spec_cA$. As $A$ is simple, $A\to A/C$ is injective, so $C=\Delta$.
The converse direction is clear, and so is (d).

For (e) let $A$ be a non-trivial sesquiad.
Bein non-trivial is equivalent to $0\ne 1$.
The Lemma of Zorn gives us a maximal congruence $C$ with $0\nsim_C 1$.
As $C$ is maximal, $A/C$ is simple.
(f) is immediate from (e).

For (g) let $A$ be simple.
Then $A_1= A\sm \{ 0\}$ is a subset of the unit group $R_A^\times$, for otherwise there would be a prime ideal $\p\subset R_A$ with $\p\cap A_1\ne\emptyset$.
Then $A\to R_A/\p$ is no longer injective, so indeed $A_1\subset R_A^\times$.
Further, as $\spec_c A$ consists of $\Delta$ alone, we get $\Ga A=A_\Delta=S_\Delta^{-1}A$.
As this is a subset of $R_A$, it follows that $S_\Delta\subset R_A^\times$, too.
In total we get $\Ga A=A_\Delta\subset R_A^\times\cup\{ 0\}$, hence $\Ga A$ is integral.
\end{proof}

\begin{definition}
We define the \e{essential spectrum} $\Ess(A)$ of a sesquiad $A$ to be the set of all $E\in\spec_cA$ such that $\Ga(A/E)$ is integral.
\end{definition}

\begin{example}
Recall the sesquiad $A=(\Z/15\Z)^\times$ of Examples \ref{Ex2.3.10}.
The sesquiad $A$ is integral, so $\Delta\in\spec_cA$, but $\Ga(A/\Delta)=\Ga A$ is not integral, so $\Delta\notin\Ess(A)$.
\end{example}

Note that if $A$ is a monoid, then 
$$
\Ess(A)=\spec_cA, 
$$
since $\Ga(A/E)=A/E$ by Proposition \ref{propcons}.

\begin{lemma}
For a given sesquiad $A$ the set $\Ess(A)$ contains all closed points of $\spec_cA$.
So $\Ess(A)$ meets every non-empty closed set.
\end{lemma}

\begin{proof}
Let $C$ be a closed point of $\spec_cA$.
Then $\Ess(A)\cap\{C\}=\Ess(A/C)$, so it suffices to show that $\Ess(A)\ne\emptyset$ for every simple $A$.
As $A$ is simple, $\Ess(A)=\spec_cA=\{\Delta\}$ by Proposition \ref{prop2.5.3}.
\end{proof}

\begin{definition}
We define the map $\sigma:\Ess(A)\to\spec_c\Ga A$ by
$$
\sigma E=\ker_c(\Ga A\to\Ga(A/E)).
$$
\end{definition}

\begin{lemma}
If we equip $\Ess(A)$ with the subspace topology of $\spec_cA$, then the map $\sigma$ is continuous.
For $E\in\Ess(A)$ one has
$$
\ga^*\sigma E=E.
$$
\end{lemma}

\begin{proof}
Let $s,t\in\Ga A$ and let $D(s,t)\subset\spec_c\Ga A$ be the corresponding open set.
Then 
$$
\sigma^{-1}(D(s,t))=\{ E\in\Ess(A): s\ne t \in\Ga(A/E)\}.
$$
Let $E,F\in\spec_cA$.
On $A_F$ the congruence $E$ induces a congruence given by
$$
\frac af\sim_E\frac bg\quad\Leftrightarrow\quad 
\exists_{h\in S_F}\ agh\sim_E bfh.
$$ 
Note that $A_F/E=0$ unless $F\in\ol E$.
In the latter case, $F$ induces a congruence on $A/E$ and one gets
$$
A_F/E\cong (A/E)_F.
$$
For everey $F\in\ol E$ we fix a representation $s(F)=\frac{a_F}{f_F}$ and $t(F)=\frac{b_F}{g_F}$.
We find
\begin{align*}
\sigma^{-1}(D(s,t))&=\{ E\in\Ess(A): \exists_{F\in\ol E}\ s(F)\nsim_E t(F)\}\\
&=\bigcup_{F\in\ol E}\bigcup_{h\in S_F}
\underbrace{D(a_Fg_Fh,b_Ff_Fh)}_{\text{open}}.
\end{align*}
So this set is open and $\sigma$ is continuous.
The last statement of the lemma is clear.
\end{proof}

We now conclude the proof of the theorem.
Let $\iota:\Ess(A)\hookrightarrow\spec_cA$ be the inclusion.
We construct a map of sheaves
$$
\sigma^\#:\CO_{\Ga A}\to \sigma_*\iota^{-1}\CO_A
$$
as follows.
Let $U\subset \spec\Ga A$ be an open set and let $V=\sigma^{-1}U$.
Let $s\in\CO_{\Ga A}(U)$.
Let $E\in\Ess(A)$.
The point evaluation map $\delta_E:\Ga A\to A_E$ localizes to $\delta_E:(\Ga A)_{\sigma_E}\to A_E$
We define the section $\sigma^\#s$ by
$$
\sigma^\#s(E)=\delta_E(s(\sigma E)).
$$
So $(\sigma,\sigma^\#)$ is a map of sesquiaded spaces $(\Ess(A),\CO_A)\to (\spec_c\Ga A,\CO_{\Ga})$.
Consider corresponding map on global sections,
$$
\Ga\sigma:\Ga\Ga A\to\iota^*\CO_A(\Ess(A)).
$$
The global sections of $\iota^*\CO_A$ are sections $\Ess(A)\to\coprod_{E\in\Ess(A)}A_E$ which are locally restrictions of sections of $\CO_A$.
But as $\Ess(A)$ meets every closed set, every global section of $\iota^*\CO_A$ extends to a unique global section of $\CO_A$, the inverse map is the restriction of global sections, so we get $\iota^*\CO_A(\Ess(A))\cong\CO_A(\spec_cA)=\Ga A$.
So $\Ga\sigma$ induces a map $\Ga\Ga A\to\Ga A$ which satisfies the condition of Proposition \ref{propcons} (c).
This means that $\Ga A$ is conservative.
\end{proof}

%: 
\section{Schemes}
\subsection{Locality}
\begin{definition}
A congruence $E$ on a sesquiad $A$ is called a \e{maximal congruence}, if it is maximal among all congruences with $0\nsim 1$.
A sesquiad $A$ is called \e{local} if it has a unique maximal congruence.
\end{definition}

\begin{lemma}
Let $E$ be a prime congruence of a sesquiad $A$.
Then $E$ induces a congruence on the localization $A_E$ which we again write as $E$ and which is given by
$$
\ker\(A_E\to S_E^{-1}R_A/S_E^{-1}I(E)\).
$$
Then we have a natural isomorphism
$$
A_E/E\cong(A/E)_\Delta.
$$
\end{lemma}

\begin{proof}
The map $A\to A_E/E$ factors over $(A/E)_\Delta$ in a unique way.
The induecd map $(A/E)_\Delta\to A_E/E$ is easily seen to be bijective as is the attached ring homomorphism $S_E^{-1}R_A/S_E^{-1}I(E)=R_{A_E/E}\to R_{(A/E)_\Delta}=S_\Delta^{-1}R_{A/E}=S_\Delta^{-1}(R_A/I(E))$.
\end{proof}

\begin{definition}
A morphism of sesquiads $\ph\colon A\to B$ is called \emph{local} if $\ph^{-1}(B^\times)=A^\times$.
Note that ``$\supset$'' always holds, so a morphism is local if and only if $\ph^{-1}(B^\times)\subset A^\times$.
\end{definition}

\begin{lemma}\label{lem2.3.1}
\begin{enumerate}[\rm (a)]
\item For a prime congruence $F$ on $B$ the localized morphism $\ph_F:A_{\ph^*F}\to B_F$ is local.
\item Let $E\in\spec_cB$ and suppose there exists $F\in\spec_cA$ and a commutative diagram
$$
\xymatrix{
A\ar[r]^\ph\ar[d]&B_E\\
A_F\ar[ur]_{\psi}
}
$$
of sesquiad homomorphisms.
If $\psi$ is local, then $F=\ph^*E$ and $\psi$ is the localization of $\ph$.
\end{enumerate}
\end{lemma}

\begin{proof}
(a) Let $a\in A_{\ph^*E}$ such that $\ph(a)\in B_E^\times$.
Then 
$$
a=\al(a_1-b_1)^{-1}\cdots(a_n-b_n)^{-1}
$$ 
with $\al,a_j,b_j\in A$ and $\ph(a_j)\nsim_E\ph(b_j)$ for every $j$, as well as $\ph(\al)\nsim_E 0$.
By definition, this means $a_j\nsim_{\ph^*E}b_j$ and $\al\nsim_{\ph^*E}0$, so that $a\in A_{\ph^*E}^\times$.

(b) The mere existence of $\psi$ implies $\ph^*E\subset F$.
Now let $a,a'\in A$ with $a\sim_F a'$ then $a-a'$ is not a unit in $A_F$, therefore $\psi(a)-\psi(a')$ is not a unit in $B_E$, so that $\psi(a)\sim_E\psi(a')$, which implies $\ph^*E= F$. The rest is clear.
\end{proof}

A \emph{sesquiaded space} is a topological space $X$ together with a sheaf $\CO_X$, called the structure sheaf, of sesquiads.
A \emph{morphism of sesquiaded spaces} $(X,\CO_X)\to (Y,\CO_Y)$ is a pair $(f,f^\#)$, where $f$ is a continuous map $f\colon X\to Y$ and $f^\#$ is a morphism of sheaves $f^\#\colon \CO_Y\to f_*\CO_X$ of sesquiads on $Y$.
Such a morphism $(f,f^\#)$ is called \emph{local}, if for each $x\in X$ the induced morphism $f_x^\# : \CO_{Y,f(x)}\to \CO_{X,x}$ is local, i.e. satisfies
$$
(f_x^\#)^* \left( E_x\right) = E_{f(x)}
$$
A \emph{isomorphism of sesquiaded spaces} is a morphism with a two-sided inverse.
An isomorphism is always local.

\subsection{Affines}
\begin{theorem}\label{thm3.2}
\begin{enumerate}[\rm (a)]
\item 
Let $A$ be a sesquiad.
Then the pair $(\spec_cA,\CO_A)$ is a sesquiaded space.
\item
Let $A,B$ be sesquiads.
If $\ph\colon A\to B$ is a morphism of sesquiads, then $\ph$ induces a local morphism of  sesquiaded spaces
$$
(f,f^\#)\colon \spec_cB\to \spec_cA,
$$
thus giving a functorial map
$$
L:\Hom(A,B)\ \to\ \Hom(\spec_cB,\spec_cA),
$$
where  on the right hand side one only admits local morphisms.
\item 
The map $L$ fits into a commutative diagram
$$
\xymatrix{
\Hom(A,B)\ar@{^(->}[r]^{L\ \ \ \ \ \ }\ar@{_(->}[dr]_{\psi} & \Hom(\spec_cB,\spec_cA)\ar@{^(->}[d]^\Ga\\
& \Hom(\Ga A,\Ga B)\ar[d]^{\cong}\\
&\Hom(\spec_c\Ga B,\spec_c\Ga A)
}
$$
where $\Ga$ is the global sections functor and $\psi$ is the map of Theorem \ref{thm2.3.8} (e), which is a bijection onto its image $\Hom_{A,B}(\Ga A,\Ga B)$.
All maps in this diagram are injective.
So, for instance, if $B$ is conservative, then all maps in the diagram are bijections.
\end{enumerate}
\end{theorem} 

\begin{example}
Before proving the theorem, we give an example showing that the map $\Ga$ above is in general not surjective.
For this let $B=(\Z/15\Z)^\times$ with addition from $\Z/15\Z$ and let $A$ be the ring $\Z/15\Z$.
Then $\Ga B=A$ and so $\Hom(\Ga A,\Ga B)=\Hom(A,A)$ contains the identity map.
On the other hand we have $\spec_cB=\{\Delta,K_3,K_5\}$, where $K_\nu$ is the congruence kernel of $B\to\Z/\nu\Z$ for $\nu=3,5$, where the generic point $\Delta$ is contained in every open set.
But $\spec_cA=\{K_3,K_5\}$ with the discrete topology. This implies that there is no continuous map from $\spec_cB$ to $\spec_cA$, hence $\Hom(\spec_cB,\spec_cA)$ is the empty set.
\end{example}

\begin{proof}
(a) is clear.
For (b) let $\ph:A\to B$ be a morphism of sesquiads.
We define $f:\spec_cB\to\spec_cA$ by $f(F)=\ph^*(F)$.
We need to define $f^\#:\CO_B\to f_*\CO_A$.

For $F\in\spec_cB$ let $\ph_F:A_{\ph^*F}\to B_F$ be the localization.
Applying this,
for any open set $U\subset\spec_cB$ we obtain a morphism
$$
f^\#:\CO_B(U)\to\CO_A(f^{-1}(U))
$$
by composing with the maps $f$ and $\ph$.
This yields a local morphism $(f,f^\#)$ of sesquiaded spaces.
We have constructed a map
$$
\eta:\Hom(A,B)\to\Hom(\spec_cB,\spec_cA).
$$
On the other hand, the global sections functor provides a map
$$
\Ga: \Hom(\spec_cB,\spec_cA)\to\Hom(\Ga A,\Ga B).
$$
By construction, the map $\Ga\circ\eta$ coincides with the map $\psi$ of Theorem \ref{thm2.3.8}, part (e).
Thus the image of $\Ga\circ\eta$ is contained in $\Hom_{A,B}(\Ga A,\Ga B)$ and $\psi^{-1}\circ\Ga\circ\eta$ is the identity on $\Hom(A,B)$, so $\psi^{-1}\circ\Ga$ is a left-inverse to $\eta$.
The injectivity of $\psi$ implies the injectivity of $L$, so 
it remains to show that $\Ga$ is injective.
For this let $(f,f^\#)$ be a local homomorphism from $\spec_cB$ to $\spec_cA$.
Then $\ph=\Ga(f^\#):\Ga A\to\Ga B$ is a morphism of sesquiads and for every $p\in\spec_cB$ there is a local morphism on the stalk $f^\#_p:\CO_{A,f(p)}\to\CO_{B,p}$.
The latter must be compatible with $\ph$ in the sense that the diagram
$$
\xymatrix{
\Ga A \ar[r]^\ph\ar[d] &\Ga B\ar[d]\\
A_{f(p)}\ar[r]^{f_p^\#}& B_p
}
$$
commutes.
Since $f^\#$ is local and the diagram commutes, the map $f_p^\#$ is the localization of the map $A\to\Ga A\tto\ph\Ga B\to B_p$ by Lemma \ref{lem2.3.1}, and so the maps $f$ and $f^\#$ are both determined by $\ph$.
The theorem is proven.
\end{proof}

\begin{example}
This is an example, where the image $\Im(\Ga)$ of the theorem is not equal to the image $\Im(\psi)$.
Let $B=\{0,1,e\}$ with $e^2=e$ and let $A=\Ga B$ (see Example \ref{Ex2.3.10}). Then $\Ga A=\Ga\Ga B=\Ga B$ and the identity $\Id:\Ga B\to\Ga B$ induces an element in $\Hom(\Ga A,\Ga B)$ which is in the image of $\Ga$, but not in the image of $\psi$.
\end{example}

\subsection{General schemes}

\begin{definition}
A sesquiaded space is called an \e{affine congruence scheme}, if it is of the form $\spec_cA$ for a sesquiad $A$.
A \e{congruence scheme} is a sesquiaded space $X$ which locally looks like  an affine, so every point $x$ has a neighborhood $U$ such that there exists a sesquiad $A$ and an isomorphism of sesquiaded spaces $U\tto\cong \spec_cA$.
\end{definition}

\begin{example}
Let $C$ be the infinite cyclic monoid generated by the element $\tau$ together with a zero, equipped with the trivial additive structure.
The congruence spectrum is the set
$$
\spec_c C=\{ \Delta,\tau\sim 0,\tau\sim 1,\tau^2\sim 1,\tau^3\sim 1,\tau^4\sim 1,\dots\},
$$
where
\begin{itemize}
\item $\tau\sim 0$ is the congruence with two equivalence classes: $\{0,\tau,\tau^2,\dots\}$ and $\{ 1\}$,
\item for $n\in\N$, we finally have the congruence $\tau^n\sim 1$ given by
$$
\tau^a\sim\tau^b\quad\Leftrightarrow\quad n\text{  divides  }a-b.
$$
The quotient monoid is isomorphic to the cyclic group $C_n$ of order $n$ extended by zero.
\end{itemize}
A non-empty subset $U$ of $\spec_cC$ is open if and only if it contains $\Delta$ and its complement $U^c$ is  a finite set which is closed under divisors, which means if $(\tau^n\sim 1)\in U^c$, then for any divisor $d$ of $n$, the element $\tau^d\sim 1$ is in $U^c$.

Take any such non-empty open set $U$.
Then two copies of $\spec_cC$ can be glued along the two copies of $U$ to give a non-affine congruence scheme.
\end{example}

\subsection{Base change to $\Z$}
The category $\ZSch$ of $\Z$-schemes is a full subcategory of the category $\CSch$ of congruence schemes.

\begin{theorem}
There is a base-change functor $(.)_\Z$ from $\CSch$ to $\ZSch$ which is the identity on $\ZSch$ and is left adjoint to the inclusion functor, so one has
$$
\Hom_\CSch(X,Y)\ \cong\ \Hom_\ZSch(X_\Z,Y)
$$
for every congruence scheme $X$ and every $\Z$-scheme $Y$.
\end{theorem}

\begin{proof}
The inclusion of categories of the category $\RINGS$ of rings into the category $\SES$ of sesquiads has a left adjoint functor given by $A\mapsto R_A$, as we have the functorial isomorphy
$$
\Hom_\SES(A,R)\ \cong\ \Hom_\RINGS(R_A,R).
$$
We define the $\Z$-base change of a sesquiad $A$ to be the ring $A_\Z=R_A$.
The $\Z$-base change of the affine congruence spectrum $\spec_cA$ is the spectrum $\spec R_A$.

Note that every affine congruence scheme $\spec_cA$, where $A$ is finitely generated as a monoid, lifts to an affine ring scheme $\spec R_A$ which comes, as an additional datum, with a natural embedding into a toric variety.
Namely, one has 
$$
\spec R_A\hookrightarrow\spec\Z_0A,
$$
and the latter was shown in \cite{F1-toric} to be a toric variety via the monoid structure of $A$.

To establish a full base change to $\Z$ we need to lift open subsets of congruence schemes.
Note that an open subset of a congruence scheme need not be a congruence scheme itself, not even contain one.
So we introduce the category $\OP_c$ of sesquiaded spaces and local morphisms, which are isomorphic to open subsets of affine congruence schemes.
Let $\OP_\Z$ be the category of $\Z$-schemes, isomorphic  to open subschemes of affine $\Z$-schemes.
We want to establish a functor 
$$
\OP_c\ \to\ \OP_\Z
$$
which extends the functor on affine schemes which is derived from the base extension functor $\SES\to\RINGS$.

Firstly, any open subset $U\subset \spec_cB$ of an affine congruence scheme can be written as
$$
U=\bigcup_{\stack{(a,b)}{D_B(a,b)\subset U}}D_B(a,b),
$$
where the union runs over all tuples $a=(a_1,\dots,a_n), b=(b_1,\dots,b_n)\in B^n$ such that the open set
$$
D_B(a,b)=D_B(a_1,b_1)\cap\dots\cap D_B(a_n,b_n)
$$
lies in $U$.
We define the $\Z$-lift of $U$ as
$$
U_\Z=\bigcup_{\stack{(a,b)}{D_B(a,b)\subset U}}D_{R_B}(a,b)\subset\spec R_B,
$$
where $D_{R_B}(a,b)$ is the set of all prime ideals $\p$ of the ring $R_B$ such that $(a_1-b_1)\cdots(a_n-b_n)\notin\p$.

Let $U\subset\spec_cB$, $V\subset\spec_cA$ be open subsets and let $(f,f^\#):(U,\CO_B)\to (V,\CO_A)$ be a local morphism of sesquiaded spaces.
In the following we define a morphism of schemes $(f_\Z,f_\Z^\#):(U_\Z,\CO_{R_B})\to(V_\Z,\CO_{R_A})$.
First we define a map $f_\Z:U_\Z\to V_\Z$.
For this let $\p\in U_\Z\subset\spec R_B$.
Let $P\in\spec_cB$ be the congruence kernel of the map $B\to R_B/\p$.
Then $P$ lies in $U$ and $f^\#$ induces a map
$$
\tilde f_P^\#: S_{A,f(P)}^{-1} R_A\to S_{B,P}^{-1}R_B.
$$
Furthermore, the localization $S_{B,P}^{-1}\p$ is a prime ideal of $S_{B,P}^{-1} R_B$, so we define
$\tilde f_\Z(\p)=(\tilde f_P^\#)^{-1}(S_{B,P}^{-1}\p)\in\spec S_{A,f(P)}^{-1} R_A$ and by $f_\Z(\p)$ we denote its inverse image in $R_A$.
By definition, it is clear that $f_\Z(\p)$ lies in $V_\Z$.
Next we show that the map $f_\Z$ is continuous.
For this let $\al\in R_A$ such that $D_{R_A}(\al)\cap V_\Z\ne\emptyset$.
We want to show that $f_\Z^{-1}(D_{R_A}(\al))$ is open in $U_\Z$.
For $\p\in U_\Z$ we have
\begin{align*}
\p\in f_\Z^{-1}(D_{R_A}(\al)\cap V_\Z) 
&\Leftrightarrow
f_\Z(\p)\in D_{R_A}(\al)\cap V_\Z\\
&\Leftrightarrow
\al\notin f_\Z(\p)\\
&\Leftrightarrow
\tilde f_P^\#(\al)\notin S_{B,P}^{-1}\p\\
&\Leftrightarrow
S_{B,P}^{-1}\p\in D_{S_{B,P}^{-1}}(\tilde f_P^\#(\al)).
\end{align*}
We infer that $f_\Z^{-1}(D_{R_A}(\al))$ is the inverse image in $\spec R_B$ under the localization map $R_B\to S_{B,P}^{-1} R_B$ of the open set $D_{S_{B,P}^{-1}R_B}(\tilde f_P^\#(\al))$, so is open indeed.
We next define a morphism of sheaves on $V_\Z$,
$$
f_\Z^\#:\CO_{R_A}\to (f_\Z)_*\CO_{R_B},
$$
by giving it on stalks.
So let $q\in V_\Z$ and let $Q\in V$ be the congruence kernel of the map $A\to R_A/\q$.
The map $f^\#$ on the stalk of $Q$ is given by
$$
f_Q^\#:A_Q=S_{A,Q}^{-1}A\ \to\ \lim_{W\supset f^{-1}(Q)}\CO_B(f^{-1}(W)).
$$
Switching to universal rings yields a map
$$
S_{A,Q}^{-1} R_A\ \to\ \lim_{W\supset f^{-1}(Q)}\CO_{R_B}(f_\Z^{-1}W_\Z).
$$
If $W$ contains $f^{-1}(Q)$, then $W_\Z$ contains $f_\Z^{-1}\q$, so we can prolong this map to
$$
S_{A,Q}^{-1}R_A\to \lim_{\tilde W\supset f_\Z^{-1}(\q)}\CO_{R_B}(f_\Z^{-1}\tilde W).
$$
As $S_{A,\q}$ contains $S_{A,Q}$, we can localize to
$$
f_{\Z,\q}^\#:\underbrace{{S_{A,\q}}^{-1} R_A}_{=\CO_{R_A,\q}}\ \to\ 
\underbrace{\lim_{\tilde W\supset f_\Z^{-1}(\q)}\CO_{R_B}(f_\Z^{-1}(\tilde W))}_{=((f_\Z)_*\CO_{R_B})_\q}
$$
This defines the map $f_\Z^\#$ on the etale-spaces of the sheaves.
It remains to show continuity, which is equivalent to showing that if $s$ is a local section of $\CO_{R_A}|_{V_\Z}$ then $f_\Z^\#\circ s$ is a local section of $(f_\Z)_*\CO_{R_B}$.
As sections are locally constant, this is easily verified.

It remains to show that $f_\Z^\#$ is local, i.e., that for $\p\in\spec R_B$ the map $f_{\Z,\p}^\#:S_{A,f_\Z^\#(\p)}R_A\to S_{B,\p}R_B$ maps non-units to non-units.
This follows from the construction of $f_\Z^\#$.

Next we define the base-change from $\CSch$ to $\ZSch$ by gluing.
Let $X$ be a congruence scheme, write $X=\bigcup_{i\in I}U_i$, where $U_i$ is in $\OP_c$.
This means that $X$ is obtained from the $U_i$ by gluing along open subsets.
The gluing maps can be lifted to corresponding open subsets of the $U_{i,\Z}$ so give gluing recipes to construct a scheme $X_\Z$.
Morphisms are treated similarly, so we get the desired base-change functor.
Since for a ring $A$ one has $R_A=A$, this functor is the identity on $\Z$-schemes.
Finally, to get the adjoint property
$$
\Hom_\CSch(X,Y)\ \cong\ \Hom_\ZSch(X_\Z,Y),
$$
it suffices to assume that $X=U\subset\spec_cB$ is an open subset of an affine.
Let $\al:U\to Y$ be a morphism and write $Y=\bigcup_iV_i$, where each $V_i=\spec A_i$ is affine.
Let $U_i=\al^{-1}(V_i)$.
By the previous, we get a morphism $\al_\Z:U_{i,\Z}\to V_{i,\Z}=\spec R_{A_i}=\spec A_i$.
These morphisms can be glued to give a morphism $\al_\Z:U_\Z\to Y$.
The map $\al\mapsto\al_\Z$ is clearly injective and it is also surjective, as the same decomposition in affines can be applied on the right hand side to start with.
\end{proof}

\subsection{Tits models}
Let $Z$ be a scheme over $\Z$.
A congruence scheme $Y$ is called a \e{model} of $Z$, if $Y_\Z\cong Z$.
Note that for a model $Y$ of $Z$ and every $\Z$-scheme $X$ one has 
$$
\Hom(X,Z)=\Hom(X,Y).
$$
Therefore, a $\Z$-scheme can, in the category of congruence schemes, be replaced by any model.

For a given conservative sesquiad $A$ and a congruence scheme $X$, we write $X(A)$ for $\Hom(\spec_cA,X)$ and call it the \e{set of $A$-valued points}.
For a ring $R$ one has
$$
X(R)= X_\Z(R).
$$
For a sesquiad $A$, the map $\ga:A\to R_A$ induces $\ga^*:\spec R_A\to\spec_cA$ and therefore a natural map
$
X(A)\ \to\ X(R_A).
$

In the nineteenfifties, Jacques Tits dreamt \cite{Tits} of a field of one element that would satisfy the formula $\GL_n(\F_1)=\Per(n)$, or more generally 
$
G(\F_1)=W_G
$
for any Chevalley group $G$, where $W_G$ is the Weyl group.
In the previous approaches to this object, \cites{Soule,F1,Haran,Durov,TV,Connes,Lorscheid},
the authors have constructed a category of schemes over $\F_1$ with base extension to $\Z$, and for given $G$ as above, a scheme $G_{\F_1}$ which base-extends to $G$ and satisfies $G_{\F_1}(\F_1)=W_G$,
where $G_{\F_1}(\F_1)$ is to be interpreted as $\Hom(\spec\F_1,G_{\F_1})$.
In our setting, $G$ and $\F_1$ can be considered as objects in the same category, so that in our setting an expression like $G(\F_1)=\Hom(\spec_c\F_1,G)$ and $G(\Z)$ likewise make sense.
However, we have the same situation as in a ring extension, where one needs to have a model of a scheme over the smaller ring.
In our case, the smaller ring is $\F_1$.
Before proceeding, we have to give a definition of $\F_1$.

\begin{lemma}
The category of sesquiads has an initial object, which we call $\F_1$.
It is the monoid $\{0,1\}$ with the trivial addition.
It satisfies 
$$
|\Hom(\F_1,A)|= 1
$$
for every sesquiad $A$.
The universal ring of $\F_1$ is $\Z$.
\end{lemma}

\begin{proof}
Clear.
\end{proof}

In \cite{F1} we defined the monoid of one element to be $\F_1$.
Recall that in \cite{F1} we worked with monoids not necessarily having a zero element, however, in the current paper we require a zero element, and we require morphisms to respect the zero.
The natural way to embed the theory of \cite{F1} into the current, is to attach a zero to every monoid.
In this way, $\F_1$ of \cite{F1} becomes the $\F_1$ of the current paper.

In the setting of congruence schemes, Tits's question becomes this:
\begin{itemize}
\item For a given Chevalley group $G_\Z$, is there a model $G$ such that the image of
$$
G(\F_1)\to G(\Z)
$$
lies in the normalizer $N(T)$ of a maximal torus $T$ such that the map
$$
G(\F_1)\to N(T)\to N(T)/T=W_G
$$
is a bijection 
$$
G(\F_1)\cong W_G?
$$
\end{itemize}
We call such a model a \e{Tits model}.
An existence proof will run as follows: in the coordinate ring $R$  of $G_\Z$ find a subsesquiad $A\subset R$ such that $R=R_A$.
Then $G=\spec_cA$ is a model for $G_\Z$.
For a suitable choice of $A$ one should find a Tits model.
We now show existence of Tits models for all series of Chevalley groups.

\begin{theorem}
The general linear group $\GL_n$, the symplectic group  $\Sp_{2n}$, and the split orthogonal $\O_n$ admit Tits models.
\end{theorem}

\begin{proof}
We start by giving a Tits model for the group $\GL_n$.
Its coordinate ring is
$$
R=\Z[X_{i,j},Y]/(\det(X)Y-1).
$$
In $R$, we consider the submonoid $A$ generated by all $X_{i,j}$, $Y$, the element $(-1)$ and the elements
$$
X_{1,1}+\dots+X_{1,n},\dots,X_{n,1}+\dots+X_{n,n}.
$$
We give $A$ the additive structure of $R$ and claim that $G=\spec_cA$ is a Tits model for $G_\Z=\GL_n$.
Note first that $\F_1$ is integral with trivial addition, hence conservative, and therefore we have 
$$
G(\F_1)=\Hom(\spec_c\F_1,\spec_cA)\ \cong\ \Hom(A,\F_1).
$$
Let $\phi:A\to\F_1$ be a sesquiad morphism.
We write $\phi(X_{i,j})=x_{i,j}$ and $\phi(Y)=y$.
Since the sum $X_{i,1}+\dots+X_{i,n}$ is defined in $A$, but $\F_1$ has trivial addition, it follows that at most one of the elements $x_{i,1},\dots,x_{i,n}$ can be non-zero.
This means that $\phi$ defines an $\F_1$-valued matrix which has at most one $1$ in every row and is zero otherwise.
Next in $R$ we have $\det(X)y=1$, so
$$
\sum_\sigma{\rm sgn}(\sigma)X_{1,\sigma(1)}\cdots X_{n,\sigma(n)}Y=1,
$$
where the sum runs over all permutations $\sigma$ in $n$ letters.
We apply $\phi$ to get
$$
\sum_\sigma\phi({\rm sgn}(\sigma))x_{1,\sigma(1)}\cdots x_{n,\sigma(n)}y=1,
$$
in particular, there exists a permutation $\sigma$ such that $x_{1,\sigma(1)}\cdots x_{n,\sigma(n)}$ is non-zero, which means that $x$ is a permutation matrix.
On the other hand, every permutation matrix defines a sesquiad morphism $\phi$ and so
$$
\Hom(A,\F_1)\cong\Per(n)= W_G.
$$
The image of $G(\F_1)$ in $G(\Z)$ is the set of permutation matrices and the claim follows.

Now for the group $\Sp_n$.
For a ring $R$, the group $\Sp(R)$ is defined to be the the group of all $2n\times 2n$ matrices $g$ over $R$ with $gJg^t=J$, where $J=\smat\ {-I}I\ $ and $I$ is the $n\times n$ unit matrix.
The group $\GL_n$ embeds as a subgroup via $g\mapsto \smat g\ \ {g^{-t}}$.
In this way, the maximal torus of diagonal matrices is mapped to a maximal torus and also the Weyl group is preserved.
The coordinate ring of $\Sp_n$ is
$$
\Z[X_{i,j}]/XJX^t-J,
$$
where $i$ and $j$ run from $1$ to $2n$ and $XJX-J$ stands for the ideal generated by all entries of this matrix.
The monoid $A$ generated by 
\begin{itemize}
\item $X_{i,j}$ for $1\le i,j\le 2n$, and
\item $X_{1,1}+\dots+X_{1,n},\dots,X_{n,1}+\dots+X_{n,n}$, as well as
\item $1+\sum_{i=1}^n\sum_{j=n+1}^{2n} X_{i,j}$ and $1+\sum_{i=n+1}^{2n}\sum_{j=1}^n X_{i,j}$
\end{itemize}
Defines a Tits model for $\Sp_n$ as is seen similar to the first example.

Finally, the group $\O_n$ is the group of all $n\times n$ matrices $g$ with $gQg^t=Q$, where 
$$
Q\ = \ \left(\begin{array}{cccc}q &  &  &   \\ & \ddots &  &   \\ &  & q &   \\  &   &   & 1\end{array}\right)
$$
and $q=\smat\ 11\ $.
The 1 in the corner only appears when $n$ is odd.
The dimension of a maximal torus equals the number $k$ of $q$'s and each element has the same shape as $Q$, except that each $q$ is replaced with a matrix of the form $\smat\ a{1/a}\ $.
From these data it is easy to construct a Tits model for this group along the lines of the previous cases.
\end{proof}

As is made clear by the comment following Problem B in \cite{Lorscheid}, it is not expected that Tits models should be group objects themselves, as that would imply a splitting of the exact sequence
$$
1\to T\to N(T)\to W_G\to 1
$$
which is defined over $\Z$. Such a splitting does not exist for instance in the case $G=\SL(2)$.

%: 
\section{Log schemes}
In this section we clarify the relation to log schemes.
To a given congruence scheme we construct a log scheme, giving a functor from the category of congruence schemes to the category of log schemes.
This functor is neither injective nor surjective on isomorphy classes.

Recall that a log scheme, as introduced 
by Kato \cite{Kato} and Fontaine/Illusie 
\cite{Illusie}, consists of a 
$\Z$-scheme $X$ together with a log structure.
\begin{itemize}
\item A \e{pre-log structure} on a $\Z$-scheme $X$ is a sheaf $\CM$ of monoids on $X$ and a morphism of sheaves of multiplicative monoids 
$\al:\CM\to \CO_X$.
\item A \e{log structure} on $X$ is a pre-log structure $\al$ such that the induced morphism $\al^{-1}(\CO_X^\times)\to\CO_X^\times$ is an ismorphism.
\end{itemize}
A scheme with a pre-log structure is also called \e{pre-log scheme}.
To a given pre-log structure $\al:\CM\to\CO_X$ there exists a universal log structure $\tilde\al:\tilde\CM\to\CO_X$ defined as the push-out of the diagram
$$
\xymatrix{
\al^{-1}(\CO_X^\times)\ar@{^(->}[r]\ar[d]_\al
&\CM\\
\CO_X^\times.
}
$$

Let $(X,\CM),(Y,\CN)$ be pre-log schemes.
A morphism $f:$ of pre-log structures is a morphism $f:X\to Y$ of schemes together with a morphism of sheaves of monoids $f^\#:f^{-1} \CN\to \CM$ such that the diagram
$$
\xymatrix{
f^{-1} \CM\ar[r]^{f^\#}\ar[d]^\al 
& \CN\ar[d]^\al\\
f^{-1}\CO_Y\ar[r]^{f^\#}
& \CO_X
}
$$
commutes.

For a given congruence scheme $X$, let $X_\Z$ be the $\Z$-base change.
Define the sheaf of monoids $\CM_X$ by
$$
\CM_X(U)=\lim_{\stack{\to}{V_\Z\supset U}}\CO_X(V),
$$
where the limit runs over all open sets $V\subset X$ such that the $\Z$-lift $V_\Z$ of $V$ contains $U$.
The natural maps $\CO_X(V)\to\CO_{X_\Z}(V_\Z)$ induce a morphism of sheaves $\al:\CM\to\CO_{X_\Z}$, i.e., a  pre-log structure.
We denote by $\Phi(X)$ the resulting pre-log-scheme.
Via the functor $\Phi$, we can consider congruence schemes as special pre-log-structures, as the following Theorem shows.

First we need some notation.
For a sesquiad $(A,R_A)$, we say that an open set $V\subset\spec R_A$ is \e{defined over} $A$, if $V$ is the base-change $U_\Z$ of some open subset $U$ of $\spec_cA$.
This is equivalent  to saying that the ideal $I(V)$ is in the image of the ideal map $I$ of Lemma \ref{lem2.1.4}.
Next assume that $V$ is defined over $A$ and let $(B,R_B)$ be another sesquiad with an open set $V'\subset R_B$ defined over $B$, say $V=U_\Z$ and $V'=U_\Z'$.
Then a morphism $\ph:V\to V'$ is \e{defined over} $(A,B)$, if $\ph$ is the lift of a morphism $U\to U'$ in the category of sesquiaded spaces.

\begin{theorem}
The functor $\Phi$ from the category $\CSch$ of congruence schemes to the category of pre-log-schemes is faithful.

Its image consist of all pre-log schemes $(Z,\CM,\al)$ such that $Z$ admits a covering $Z=\bigcup_iV_i$ by affines $V_i=\spec R_i$ such that $\al:\CM(V_i)\to R_i$ is injective, its image generates the ring $R_i$, each open set $V_{i,j}\subset \spec R_i$ given by the gluing, is defined over $\al(\CM(V_i))\cong\CM(V_i)$ and each gluing morphism $V_{i,j}\to V_{j,i}$ is defined over $(\CM(V_i),\CM(V_j))$.
\end{theorem}

\begin{proof}
Let $\PlogSch$ denote the category of pre-log-schemes.
We have to show that $\Phi:\Hom_{\CSch}(X,Y)\to \Hom_{\PlogSch}(\Phi X,\Phi Y)$ is injective.
So suppose $\Phi(f)=\Phi(g)$.
Then the lifts $f_\Z$ and $g_\Z$ agree and so the continuous maps $f$ and $g$ on the space $X$ agree.
Since also the induced maps on $f_\Z^{-1}\CM_Y\to\CM_X$ agree, we have $f^\sharp=g^\sharp$, therefore $\Phi$ is faithful.

The image of $\Phi$ is contained in the denominated class, as any affine covering of $X$ yields a covering of the given type of $\Phi X$.
For the converse note that the gluing data for $Z$ as in the theorem induces gluing data for the sesquiads $(\al(\CM(V_i)),R_i)$, given a preimage of $Z$ under $\Phi$.
\end{proof}

\begin{example}
We give an example to show that this pre-log structure is not a log structure in general.
Consider the sesquiad $A=\{ 0,1\}\cup 2^\N$ consisting of all powers of $2$ inside $\Z$, which is its universal ring.
Let $X$ be the congruence scheme given by $\spec_cA$ and let $X_\Z$ be its base change.
As $\Z$ is the universal ring of $A$, we have $X_\Z=\spec\Z$.
Let $\CM$ denote the corresponding log structure on $X_\Z=\spec\Z$.
Let $\eta$ denote its generic point, then $\CM_\eta$ equals the localization of $X_2$ by the diagonal, a monoid which is a submonoid of $S^{-1}X_2\subset\Q$, where $S=\Z\sm\{ 0\}$.
As $S^{-1} X_2$ consists of all rational numbers with powers of two as enumerators, its unit group equals $2^\Z$, and this is a proper subset of $\Q^\times=\CO_{X_\Z,\eta}^\times$.
So the map $\CM_\eta^\times\to\CO_{X_\Z,\eta}^\times$ is not surjective and thus the pre-log structure is not a log structure.
\end{example}

%: 
\section{Berkovich subdomains}
\subsection{Tame sets}

\begin{definition}
Let $(A_i)_{i\in I}$ be a directed family of sesquiads with universal rings $R_i$.
For $i\le j$, let $\ph^j_i:R_i\to R_j$ be the corresponding ring homomorphism.
Then $A=\lim_\to A_i$ is a  sesquiad with universal ring $R_A=\lim_\to R_i$.
\end{definition}

Let $(X,\CO)$ be a congruence scheme.
For an arbitrary subset $S\subset X$ we define the sesquiad
$$
\CO(S)\df\lim_{\stack{\to}{U\supset S}}\CO(U),
$$
where the limit extends over all open subsets of $X$ containing $S$.
It $S\subset T\subset X$, then there is a natural restriction morphism
$\CO(T)\to\CO(S)$.

\begin{definition}
Let $A$ be a sesquiad.
A set $\emptyset\ne T\subset\spec_cA$ is called \e{tame}, if the image of the map $\spec_c\CO(T)\to\spec_cA$ lies in $T$.

An open set $U\ne \emptyset$ is called \e{basic}, if it is of the form $U=D(a,b)$ for some $a,b\in A^n$.
\end{definition}

\begin{example}
We give an example of a non-tame set.
Let $A$ be the monoid $\{0,1,\tau\}$ with $ \tau^2=1$.
Then $\spec_cA=\{\Delta,c\}$ where the closed point $c$ is given by $\tau\sim_c1$.
The only open set that contains $c$, is $U=\spec_cA$.
Therefore the set $S=\{ c\}$ is not tame, as $\CO(S)=A$.
\end{example}

\begin{proposition}\label{proptame}
Let $A$ be a sesquiad and $U\subset\spec_cA$ be an open set.
\begin{enumerate}[\rm (a)]
\item \label{proptame0} Any intersection of tame sets is tame.
\item \label{proptame1}If $A$ is  a ring, then every open set $U\subset \spec_cA$ is tame.
\item \label{proptame2}If for every proper ideal $I$ of $R_{\CO(U)}$ the ideal $I\CR(U)$ is proper, then $U$ is tame.
\item \label{proptame3}If the ring extension 
$\CR(U)/R_{\CO(U)}$ is integral, then $U$ is tame.
\item \label{proptame4}If $U$ is basic, then $U$ is tame.
\item \label{proptame5}If $U$ is an affine subscheme of $\spec_cA$, then $U$ is tame.
\item \label{proptame6}If $U$ has finitely many connected components, each of which is tame, then $U$ is tame. 
\end{enumerate}
\end{proposition}

\begin{proof}
(\ref{proptame0}) Let $T_i$ be tame for each $i\in I$ and let $T=\bigcap_{i\in I}T_i$.
The map $A\to\CO(T)$ factors over $\CO(T_i)$, therefore, as $T_i$ is tame, the image of $\spec_c\CO(T)$ lies in $T_i$.
Since this holds for every $i\in I$, it follows that $T$ is tame.

(\ref{proptame1}) is classical.
For (\ref{proptame2}) note that the condition is equivalent to saying that the map $\Mspec \CR(U)\to\Mspec R_{\CO(U)}$ is surjective.
By Lemma \ref{meetsevery} it follows that the map $\spec\CO_{R_A}(V)=\spec \CR(U)\to\spec_c\CO(U)$ meets every non-empty closed set. Let $\rho:A\to R_A$ the inclusion and $V=(\rho^*)^{-1}(U)$.
By Lemma \ref{lem2.3.10} we have $\CR(U)\cong\CO_{R_A}(V)$.
As the map $\CO(U)\to\CR(U)=\CO_{R_A}(V)$ factors over $R_{\CO(U)}$, the map $\spec\CO_{R_A})(V)\to\spec_c\CO(U)$ factors over $\spec R_{\CO(U)}$.
The tameness of $V$ implies that the following diagram
$$
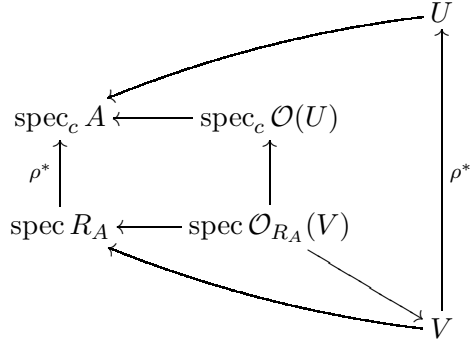
\xymatrix{
&&U\ar@/_/[lld]\\
\spec_cA & \spec_c\CO(U)\ar[l]\\
\spec R_A\ar[u]^{\rho^*} & \spec\CO_{R_A}(V)\ar[l]\ar[u]\ar[dr]\\
&&V\ar@/^/[llu]\ar[uuu]_{\rho^*}
}
$$
commutes.
The commutativity of the diagram and the fact that $V$ is tame by (a) implies that $\spec\CO_{R_A}(V)$ maps into $U$.
Let now $C\subset\spec_c\CO(U)$ be the inverse image of $\spec_cA\sm U$, then $C$ is closed.
So if $C\ne\emptyset$, then there is $x\in\spec\CO_{R_A}(V)$ mapping into $C$, which contradicts the fact that $\spec\CO_{R_A}(V)$ maps into $U$, so $C$ is empty, i.e., $U$ is tame.

(\ref{proptame3}) Let $S/R$ be an integral ring extension and let $I$ be a proper ideal of $R$.
We claim that $IS$ is a proper ideal of $S$.
It suffices to assume that $I$ is a maximal ideal.
If $IS$ was not proper, then $1\in IS$, so $1=\sum_{j=1}^nx_js_j$ with $x_j\in I$ and $s_j\in S$.
Then $s_1,\dots,s_n$ generate a finite ring extension, so it suffices to assume $S/R$ is finite.
As $R\to S$ is injective, the map $\spec S\to\spec R$ is dominant.
As $S/R$ is finite, this map is closed.
Together it follows that $\spec S\to\spec R$ is surjective, so there is an ideal $\p$ of $S$ with $\p\cap R=I$.
But then $\p\supset IS$, so $IS$ is proper.

(\ref{proptame4}) Let $U=D(a,b)$ be basic.
Let $f=(a_1-b_1)\cdots(a_n-b_n)$ and let $S_f$ be the submonoid of $R_A$ generated by $f$.
Then 
\begin{align*}
\CR(U)&=\CO_{R_A}((\rho^*)^{-1}(U))= \CO_{R_A}((\rho^*)^{-1}(D(a,b)))\\
&= \CO_{R_A}(D(f))= S_f^{-1}R_A= R_{S_f^{-1}A}.
\end{align*}
There is a natural map $S_f^{-1}A\to\CO(U)$ which induces a map $\CR(U)=R_{S_f^{-1}A}\to R_{\CO(U)}$.
The latter is inverse to the embedding $R_{\CO(U)}\hookrightarrow \CR(U)$.
So the map $R_{\CO(U)}\to\CR(U)$ is an isomorphism, so by (c) the set $U$ is tame.

(\ref{proptame5}) $U$ being an affine subscheme means that there is a sesquiad morphism $\ph:A\to B$ such that $\ph^*:\spec_cB\to\spec_cA$ is a homeomorphism $\spec_cB\tto\cong U$ and the corresponding sheaf homomorphism $\CO_A|_U\to \CO_B$ is an isomorphism.
We have a commutative diagram
$$
\xymatrix{
A\ar[r] \ar[d] & B\ar[d]\\
\Ga A=\CO(\spec_cA)\ar[r] & \Ga B=\CO(U).
}
$$
It follows that the map $A\to\CO(U)$ factors over $\ph$, which implies  that the map $\spec_c\CO(U)\to\spec_cA$ factors over $\spec_cB=U$, which means that $U$ is tame.

(g) Assume $U=U_1\cup\dots\cup U_n$ is the decomposition into connected components.
A connected component is always a closed set.
Therefore $U_2\cup\dots\cup U_n$ is closed, so $U_1$ is open as are all components.
It follows that $\CO(U)=\CO(U_1)\times\dots\times\CO(U_n)$ so the spectrum of $\CO(U)$ is the disjoint union of the spectra of the $\CO(U_j)$.
The claim follows.
\end{proof}

\subsection{Monoids}
We consider a monoid as a sesquiad with trivial addition.
For a sesquiad $A$, an open set $U\subset\spec_cA$ is called \e{monoidal}, if the sesquiad $\CO(U)$ is a monoid.

\begin{definition}[Berkovich]
Let $A$ be a monoid.
An \e{affine subdomain} is a subset $U\subset\spec_cA$ together with a morphism $\ph_U:A\to A_U$ of monoids such that 
$$
\ph_U^*(\spec_cA_U)\subset U,
$$ 
and the following universal property holds: For every monoid morphism $\ph:A\to B$ with $\ph^*(\spec_cB)\subset U$ there exists a unique morphism of monoids $A_U\to B$ making the diagram
$$
\xymatrix{
A\ar[r]^{\ph_U}\ar[dr]_{\ph} & A_U\ar@{.>}[d]^{\exists !}\\
& B
}
$$
commutative.
\end{definition}

\begin{conjecture}
Let $A$ be a monoid. An open subset $U\subset \spec_cA$ is an affine subdomain if and only if $\CO(U)$ is a monoid, i.e., has trivial addition.
In this case one has $A_U=\CO(U)$.
\end{conjecture}

\begin{definition}
We call a sesquiad $A$ a \e{division sesquiad}, if $A\sm\{0\}$ is a group, i.e., if $A=A^\times\cup\{0\}$.

For an integral sesquiad $A$, let $S_A$ be the submonoid $A\sm\{ 0\}$ and define the \e{quotient sesquiad} of $A$ as
$$
Q(A)\df S_A^{-1}A.
$$
This then is a division sesquiad as $Q(A)=G\cup\{0\}$, where $G$ is a group and $A$ injects into $Q(A)$.

For an arbitrary sesquiad $A$, let $\CZ$ denote the zero class map from $\spec_cA$ to $\spec_zA$.
For $E\in\spec_cA$ let
\begin{align*}
\what E&= \Im\(\spec_c Q(A/E)\to\spec_cA\)\\
&= \ol E\cap\CZ^{-1}(\CZ(E)).
\end{align*}
A set $C\subset\spec_cA$ is called \e{semi-closed}, if $E\in C\Rightarrow\what E\subset C$.
For a set $U\subset\spec_cA$, define the \e{semi-closed core} as
$$
SC(U)\df\{E\in U:\what E\subset U\}.
$$
\end{definition}

\begin{lemma}
Let $A$ be a monoid and $U\subset\spec_cA$ be open.
Then a prime congruence $E$ is in $SC(U)$ if and only if the whole $\CZ$-fibre $\CZ^{-1}(\CZ(E))$ lies in $U$. 
\end{lemma}

\begin{proof}
If the $\CZ$-fibre of $E$ lies in $U$, then $\what E$ lies in $U$, so $E\in SC(U)$.
The other way round, let $E\in SC(U)$.
The unique closed point $C\in\what E$ lies in $U$, so there are $a,b\in A^n$ such that $C\in D(a,b)\subset U$.
As $C$ has only two equivalence classes,  we can assume that $a_j\sim_C1$ and $b_j\sim_C 0$ for every $j=1,\dots,n$.
A prime congruence $F$ lies in the $\CZ$-fibre of $E$ if  and only of $[0]_F=[0]_C$, therefore, any such $F$ will be in $D(a,b)$.
\end{proof}

\begin{proposition}
Let $A$ be a monoid.
\begin{enumerate}[\rm (a)]
\item Let $A$ be integral and let $\p\in\spec_zA$.
The fibre $\CZ^{-1}(\p)$ is in bijection with the set of all subgroups $H$ of the quotient sesquiad $Q(A/\p)$ as follows.
For given $H$, define a prime congruence $E_H$ by
$$
x\sim 0\ \Leftrightarrow\ z\in\p,\quad\text{and otherwise}\quad x\sim y\ \Leftrightarrow\ x^{-1} y\in H.
$$
Then the map $H\mapsto E_H$ is the claimed bijection.
\item Let $\ph:A\to B$ be morphism of monoids.
Then $\Im\ph^*$ is semi-closed.
\item
Let $U\subset\spec_cA$ be an affine subdomain with structure morphism $\ph_U:A\to A_U$.
Then we have
$$
\Im(\ph_U^*)=SC(U).
$$
\item Let $A$ be a monoid and let $U\subset \spec_cA$.
Then $U$ is an affine subdomain if and only if $SC(U)$ is an affine subdomain.
In this case one has $\ph_U=\ph_{SC(U)}$.
\item Let $A$ be a monoid and $U$ a semi-closed subset of $\spec_cA$.
If $U$ is an affine subdomain, then $U$ is tame and we have $\CO(U)\cong\Ga A_U$.
Further, the map $\ph_U^*:\spec_cA_U\to U$ is a bijection.
\end{enumerate}
\end{proposition}

\begin{proof}
(a) Let $\ph:\{H\subset Q(A/\p)\}\to\CZ^{-1}(\p)$ the ensuing map.
We show injectivity of $\ph$: Let $H\ne H'$, say $h\in H\sm H'$.
Then there exists $x,y\in A\sm\p$ with $h=x^{-1} y$ so that $x\sim_{\ph(H)}y$, but $x\nsim_{\ph(H')}y$, which shows $\ph(H)\ne\ph(H')$.
For surjectivity, Let $E\in\CZ^{-1}(\p)$ and let $H$ be the set of all $x^{-1}y\in Q(A/ \p)^\times$ where $x\sim_Ey$ in $A\sm\p$.
Then $H$ is a group and $\ph(H)=E$ by definition.

(b) Let $E\in\Im\ph^*$, say $E=\ph^*F$.
We get a commutative diagram
$$
\xymatrix{
A\ar[r]\ar[dr] & B/F\ar[r] & Q(B/F)\\
& A/E\ar@{^(->}[u]\ar[r] & Q(A/E).\ar@{^(->}[u]
}
$$
It follows from (a) that the map $\spec_cQ(B/F)\to\spec_cQ(A/E)$ is surjective.
This implies the claim.

(c) If $E\in\Im(\ph_U^*)$, then $\what E\subset \Im(\ph_U^*)\subset U$ by part (b).
For the converse direction, let $E\in U$ with $\what E\subset U$.
Let $\ph:A\to A/E$, then $\ph^*(\spec(A/E))\subset U$, therefore $\ph$ factors over $\ph_U$, so $E\in\Im(\ph_U^*)$.

(d) If $U$ is an affine subdomain, then $\Im(\ph_U)\subset SC(U)$ by (d).
Therefore $SC(U)$ is an affine subdomain with $\ph_{SC(U)}=\ph_U$.
For the converse, assume that $SC(U)$ is an affine subdomain and let $\ph:A\to B$ with $\Im(\ph^*)\subset U$.
By (b), the image of $\ph^*$ lies in $SC(U)$, therefore $\ph$ factors uniquely over $\ph_{SC(U)}$, which means that $U$ is an affine subdomain with $\ph_U=\ph_{SC(U)}$.

(e) We show that if $U$ is a semi-closed affine subdomain, then it is monoidal and tame.
We first show that the map $\ph_U^*:\spec_cA_U\to U$ is a bijection.
Let $E\in U$, then the image of $\spec_cQ(A/E)$ lies in $U$, so the map $\spec_cQ(A/E)\to U$ factors over $\spec_cA_U$, hence $E$ lies in the image of $\spec_cA_U\to U$, which therefore is surjective.
For the injectivity let $E=\ph_U^*F_1=\ph_U^*F_2$.
We get the commutative diagram with solid arrows,
$$
\xymatrix{
A\ar[r]\ar[d] & Q(A/E)\ar[d]\\
A_U\ar@{.>}[ru]\ar[r] & Q(A_U/F_1).
}
$$
As $\spec_cQ(A/E)$ maps into $U$, there is a dotted arrow making the upper triangle commutative.
The uniqueness part in the universal property of $A_U$ implies that also the lower triangle commutes.
Therefore, the right vertical map is surjective and $F_1$ is uniquely determined by $E$, so $F_1=F_2$.

As $\spec_cA_U$ maps into $U$, the map $\ph_U$ induces a continuous map $f=\ph_U^*:\spec_cA_U\to U$, which we have shown to be a bijection, and a map of sheaves $f^\#:\CO_U\to f_*\CO_{A_U}$.
We show the latter to be an isomorphism.
It is sufficient to show this on stalks.
So let $E\in U$ and let $F\in\spec_cA_U$ the unique element with $\ph_U^*F=E$.
Taking global sections, we get a map from $\CO(U)$ to $\Ga A_U$ which fits into the following commutative diagram
$$
\xymatrix{
A\ar[rr]\ar[dr]\ar@/_/[ddr] && A_U\ar[d]\\
& \CO(U)\ar[r]\ar[d] & \Ga A_U\ar[d]\\
& A_E\ar[r] & A_{U,F}.
}
$$
The image of $\spec_cA_E\to\spec_cA$ consists of the intersection of all $D(a,b)$ which contain $E$, so it is the intersection of all open neighborhoods of $E$.
This means that for each $E\in U$, the image of $\spec_cA_E\to\spec_cA$ is also contained in $U$, which in turn implies that the map $A\mapsto A_E$ factors over $A_U$.
We get this diagram
$$
\xymatrix{
A\ar[rr]\ar[dr]\ar@/_/[ddr] && A_U\ar[d]\ar@{.>}[ddl]\\
& \CO(U)\ar[r]\ar[d] & \Ga A_U\ar[d]\\
& A_E\ar[r] & A_{U,F}.
}
$$
The dotted arrow localizes to a map $A_{U,F}\to A_E$, which, by the commutativity of the diagram, inverts to lower horizontal arrow, so the stalks are isomorphic.
Taking global sections, we find that the middle horizontal map in the diagram is an isomorphism.
Therefore the map from $A$ to $\CO(U)$ factors over $A_U$, so the map $\spec_c\CO(U)\to\spec_cA$ factors over $\spec_cA_U$, hence $U$ is tame.
\end{proof}

\begin{examples}
\item We give an example of an open set $U$, such that the semi-closed core $SC(U)$ is no longer open.
This is the free monoid $C$ with one generator $\tau$.
Let $U=D(1,\tau)$, then $SC(U)$ consists of one closed point $\{\tau\sim 0\}$ only.
The set $U$ is even an affine subdomain.
In this case one has $A_U=\{ 0,1\}$ and the map $C\to A_U$ sends $\tau$ to zero.
Note that $A_U\ne\CO(U)=C[(1-\tau)^{-1}]$ and that $\CO(U)$ is not a monoid.
\item The empty set $U=\emptyset$ is an open affine subdomain.
We have $A_U=\{1\}$ in this case.
Here $1=0$ and so $A_U$ has empty spectrum.
\item Let $A=\{ 0,1,a,b,ab\}$ with $a^2=b^2=1$.
Then $\spec_cA$ is in bijection with the subgroups of $A^\times$.
Let $U=\{\Delta\}$, then $U$ is an open affine subdomain with $A_U=\{1\}$ again, since the only monoid $B$ for which there exists a monoid morphism $A\to B$ which maps $\spec_cB$ into $U$, is the trivial monoid $B=\{1\}$ which has empty spectrum.
\end{examples}

%: 
\section{Zeta functions}
Recall the Hasse-Weil zeta function of a scheme $Z$ of finite type over $\Z$ is, for $\Re(s)>>0$, defined as
$$
\zeta_Z(s)=\prod_{z\in|Z|}\frac1{1-N(z)^{-s}},
$$ 
where $|Z|$ is the set of closed points of $Z$ and $N(z)$ is the cardinality of the residue field at $z$.

\subsection{Closed points}
In order to transfer the notion of zeta function to congruence schemes, we will first investigate the set of closed points.

\begin{examples}
\item Let $F_\tau$ be the sesquiad $\{0,1,\tau,\tau^2,\dots\}$ freely generated by $\tau$ with the trivial addition.
Recall that $\spec_c F_\tau=\{ \Delta, \tau\sim 0,\tau\sim 1,\tau^2\sim 1,\dots\}$.
The closed points are $\tau\sim 0$ and $\tau\sim 1$.
\item For a prime number $p$ let $X_p$ denote the sesquiad $\{0,1,p,p^2,\dots\}$ inside the ring $\Z$.
For a natural number $m\ge 2$, let $C_p(m)$ denote the congruence
$$
C_p(m)=\ker(X_p\to\Z/m).
$$
Any congruence $\ne\Delta$ of $X_p$ is of this form.
If $p$ is not invertible modulo $m$, then the image of $X_p$ in $\Z/m$ can only be integral, if $m=p$.
So we get
$$
\spec_cX_p=\{\Delta,C_p(p)\}\cup\{C_p(m)\in\N:m>1\text{ and }(m,p)=1\}.
$$
As $(\Z/m)^\times$ is a finite group, there exists a smallest $n\in\N$ with $p^n\equiv 1\mod(m)$, so $m$ divides $p^n-1$.
It follows that the map $\Z\to\Z/(p^n-1)$ induces the same congruence on $X_p$ as $\Z\to\Z/m$, therefore we can identify $\spec_cX_p$ with the set $\{\Delta,C_p(p)\}\cup\{n=1,2,\dots\}$ via $n\mapsto C_p(p^n-1)$.
Let $1<d<n$ be a divisor of $n$, then by
$$
(p^d-1)(1+p^d+p^{2d}+\dots+p^{d(\frac nd-1)})=p^n-1
$$
we infer that $p^d-1$ divides $p^n-1$, so we get a ring homomorphism $\Z/(p^n-1)\to\Z/(p^d-1)$ and thus a sesquiad homomorphism $X_p/C_n\to X_p/C_d$ which means that $C_n$ is not a closed point.
We conclude that if $C_n$ is a closed point, then either $n$ is a prime  or $n=1$ and $p\ne 2$.
Therefore we have
$$
|X_p|=\begin{cases}\{C_p(p)\}\cup\{ C_p(p^l-1):l\text{ prime or }l=1\}& p\ne 2,\\
\{C_2(2)\}\cup \{ C_2(2^l-1):l\text{ prime}\}&p=2.\end{cases}
$$
\end{examples}

For a congruence scheme $X$ and a point $x\in X$,
we have a natural congruence $E_x$ on the local sesquiad $\CO_{X,x}$ induced by $x$.
We define the \e{residue sesquiad} of the point $x$ as
$$
\res(x)=\CO_{X,x}/E_x.
$$
If $X=\spec_cA$, one has
$$
\res(x)\cong\(A/E_x\)_\Delta.
$$
\begin{definition}
A point $x$ in a congruence scheme $X$ is called a \e{$\Z$-point}, if the ring $R_{\res(x)}$ is integral.
\end{definition}

\begin{lemma}
\begin{enumerate}[\rm (a)]
\item Let $A$ be an integral sesquiad.
Then one has
$$
R_A\text{ integral}\quad\Leftrightarrow\quad R_{A_\Delta}\text{ integral}.
$$
\item Let $A$ be a sesquiad and $x\in X=\spec_cA$.
Then $x$ is a $\Z$-point if and only if $R_{A/E_x}$ is integral.
\end{enumerate}
\end{lemma}

\begin{proof}
(a)
One has
$$
R_{A_\Delta}= R_{S_\Delta^{-1}A[\Delta]}= S_\Delta^{-1}R_{A[\Delta]}=S_\Delta^{-1}R_A.
$$
If $R_A$ is integral, then so is $S_\Delta^{-1} R_A$.
For the converse direction, assume $S_\Delta^{-1} R_A$ is integral and $\al\beta=0$ in $R_A$.
Then  in the ring $S_\Delta^{-1}R_A$ one has $\al=0$ or $\beta=0$.
Let's assume $\al=0$ holds in this localization, then there exists $s\in S_\Delta$ such that $\al s=0$ holds in $R_A$.
This means that there exist $a,b\in A^n$ for some $n$ such that $\Delta\in D(a,b)$, so $D(a,b)\ne\emptyset$, and
$$
\al(a_1-b_1)\cdots (a_n-b_n)=0.
$$
This means that $\al$ is the zero section on $D(a,b)$.
As $A$ is integral, $X=\spec_cA$ is irreducible, hence the open set $D(a,b)$ is dense, so $\al=0$.

(b) We have $\res(x)=\CO_{X,x}/E_x=A_{E_x}/E_x=(A/E_x)_\Delta$. So the claim follows from part (a).
\end{proof}

\begin{definition}
A point $x$ of a congruence scheme $X$ is called \e{$\Z$-closed}, if $R_{\res(x)}$ is a field.
In this case we call
$$
\kappa(x)= R_{\res(x)}
$$
the \e{residue field} of the point $x$.
\end{definition}

\begin{example}
Condsider $X_2$.
The closed point $\tau\sim 0$ is also $\Z$-closed, since its universal ring is $\F_2$.
The closed point $\tau^p\sim 1$ is $\Z$-closed if and only if $\Z/(2^p-1)$ is a field and this is equivalent to $2^p-1$ being a Mersenne-prime.
\end{example}

Let $X$ be a congruence scheme and denote by $|X|_\Z$ the set of $\Z$-closed points.
For $x\in |X|_\Z$ we set
$$
N(x)=\begin{cases}|\kappa(x)|&\text{if }\kappa(x)\text{ is finite,}\\
\infty&\text{otherwise.}\end{cases}
$$
We define the \e{Hasse-Weil-zeta-function} of $X$ by 
$$
\zeta_X(s)=\prod_{x\in|X|_{\Z}}\frac1{1-N(x)^{-s}},
$$
if the product converges for some $s\in\C$, where we interpret $N(x)^{-s}$ to be zero if $N(x)=\infty$.

\begin{example}
We have
\begin{align*}
\zeta_{X_2}(s)&=\frac1{1-2^{-s}}\prod_{p:2^p-1\text{ prime}}\frac1{1-(2^p-1)^{-s}}\\
&=
\frac1{1-2^{-s}}\prod_{q\text{ Mersenne prime}}\frac1{1-q^{-s}}.
\end{align*}
\end{example}

%: 
\section{Closing remarks}
\subsection{Relations to blueprints and blue schemes}\label{relation}
In this paper, sesquiads have additions from rings, but, as mentioned in Section  \ref{alternative}, one can do the same with semi-rings leading to \e{semi-sesquiads}.
At the moment, it is not clear, whether the theory of congruence schemes can be extended to semi-sesquiads. 
The construction of a structure sheaf doesn't extend straightforwardly.

In \cite{blueprints}, Oliver Lorscheid has defined the notion of a blueprint and a blue scheme.
A blueprint is the same as a semi-sesquiad and a blue scheme is what one might call a Zariski-scheme over semi-sesquiads.
Here spectra are built from ideals and not congruences.
Let's call a blue scheme \e{special}, if it has a covering of affines which are Zariski-spectra of sesquiads.
A sesquiad defines two objects: an affine Zariski scheme and an affine congruence scheme.
The zero class map $E\mapsto [0]_E$ gives a fibration of the congruence scheme over the Zariski scheme.
Gluing Zariski schemes is compatible with this fibration, so every special blue scheme extends to a congruence scheme, but not every congruence scheme is obtained in this way, as in the congruence topology there are more possibilities for gluing.
The category of special blue schemes becomes a subcategory of the category of congruence schemes in this way.
It is an interesting question, whether it is a full subcategory or not.

\subsection{Open questions}
\begin{definition}
We call a congruence scheme $X$ \e{primal} if it has an open covering $X=\bigcup_{i\in I}U_i$ of affine open subschemes such that for each $i\in I$ we have  $U_i=\spec_cA_i$ for a monoid $A_i$.
\end{definition}
Let $A,B$ be   sesquiads.
\begin{enumerate}
\item Is the map $\ga^*:\spec_c\Ga A\to\spec_cA$ injective?
\item Is every open subset $U\subset\spec_cA$ tame?
This would imply that every set $S\subset\spec_cA$ with
$$
S=\bigcap_{U\supset S}U
$$
is tame, which is the best one can hope for.
\item Is there a description as in \cite{F1-toric} of the category of all $\Z$-lifts of  primal congruence schemes?
\item Is it true that $U\subset \spec_cA$ is an affine subdomain if and only if it is an affine subdomain of $\spec_c(A^\red)$?
Is then $(A_U)^\red=(A^\red)_U$?
\item Is $A^\red\to\prod_{E\in\spec_cA}A/E$ an embedding?
\item Is the ring $R_{A^\red}$ reduced?
\item Let $A$ be a monoid such that $\spec_cA$ is connected.
Is $A$ conservative?
\end{enumerate}

{\small Mathematisches Institut, Auf der Morgenstelle 10, 72076 T\"ubingen, Germany\\
\tt deitmar@uni-tuebingen.de}

\today

\end{document}